\documentclass[article, 10pt, oneside, reqno]{amsart}

\usepackage[utf8]{inputenc}
\usepackage{amsfonts,amssymb,amsthm,amsmath}
\usepackage{mathrsfs}
\usepackage{graphicx,graphics}
\usepackage{color}
\usepackage{tikz}
\usepackage{tikz-cd}
\usetikzlibrary{nfold}
\usepackage{comment}
\usepackage{subfiles}
\usepackage{stmaryrd}
\usepackage{hyperref}
\usepackage{todonotes}
\usepackage{cleveref}

\hypersetup{colorlinks = true, linkbordercolor = {white}, linkcolor = {blue}, citecolor = {blue}}

\theoremstyle{plain}
\newtheorem{thm}{Theorem}[subsection]
\newtheorem{lem}[thm]{Lemma}
\newtheorem{prop}[thm]{Proposition}
\newtheorem{corol}[thm]{Corollary}
\newtheorem{conj}[thm]{Conjecture}
\theoremstyle{definition}
\newtheorem{defi}[thm]{Definition}

\theoremstyle{remark}
\newtheorem{rque}[thm]{Remark}
\newtheorem{rques}[thm]{Remarks}

\newcommand{\F}{\mathbb{F}}

\newcommand{\Fq}{\F_q}

\newcommand{\Fqb}{\overline{\F}_q}
\newcommand{\Flb}{\overline{\F}_{\ell}}

\newcommand{\Z}{\mathbb{Z}}

\newcommand{\Zlb}{\overline{\Z}_{\ell}}

\newcommand{\Q}{\mathbb{Q}}

\newcommand{\Qlb}{\overline{\Q}_{\ell}}

\newcommand{\RR}{\mathrm{R}}

\newcommand{\Gm}{\mathbb{G}_m}

\newcommand{\Gbf}{\mathbf{G}}
\newcommand{\Bbf}{\mathbf{B}}
\newcommand{\Ubf}{\mathbf{U}}
\newcommand{\Tbf}{\mathbf{T}}
\newcommand{\Wbf}{\mathbf{W}}
\newcommand{\Nbf}{\mathbf{N}}

\newcommand{\Pbf}{\mathbf{P}}
\newcommand{\Vbf}{\mathbf{V}}
\newcommand{\Lbf}{\mathbf{L}}

\newcommand{\Mbf}{\mathbf{M}}
\newcommand{\Qbf}{\mathbf{Q}}
\newcommand{\Rbf}{\mathbf{R}}

\newcommand{\UP}{\Ubf_{\Pbf}}
\newcommand{\UQ}{\Ubf_{\Qbf}}

\newcommand{\pfrak}{\mathfrak{p}}

\newcommand{\Hom}{\mathrm{Hom}}

\newcommand{\id}{\mathrm{id}}

\newcommand{\Tr}{\mathrm{Tr}}

\newcommand{\End}{\mathrm{End}}
\newcommand{\Mod}{\mathrm{Mod}}

\newcommand{\Frob}{\mathrm{F}}

\newcommand{\Ad}{\mathrm{Ad}}

\newcommand{\opp}{\mathrm{op}}

\newcommand{\Br}{\mathrm{Br}}

\newcommand{\pr}{\mathrm{pr}}
\newcommand{\RGamma}{\mathrm{R\Gamma}}

\newcommand{\DD}{\mathrm{D}}
\newcommand{\pt}{\mathrm{pt}}

\newcommand{\Stab}{\mathrm{Stab}}

\newcommand{\Rep}{\mathrm{Rep}}

\newcommand{\rev}{\mathrm{rev}}

\newcommand{\Ccal}{\mathcal{C}}

\newcommand{\Ocal}{\mathcal{O}}

\newcommand{\Reg}{\mathrm{Reg}}

\newcommand{\Alg}{\mathrm{Alg}}

\newcommand{\Hcal}{\mathcal{H}}

\newcommand{\Lcal}{\mathcal{L}}

\newcommand{\Mcal}{\mathcal{M}}

\newcommand{\tr}{\mathrm{tr}}
\newcommand{\ch}{\mathrm{ch}}

\newcommand{\Xcal}{\mathcal{X}}

\newcommand{\hc}{\mathrm{hc}}
\newcommand{\Ch}{\mathrm{Ch}}

\newcommand{\geo}{\mathrm{geo}}

\newcommand{\Int}{\mathrm{Int}}

\newcommand{\Acal}{\mathcal{A}}
\newcommand{\Ncal}{\mathcal{N}}

\newcommand{\unit}{\mathrm{unit}}
\newcommand{\counit}{\mathrm{counit}}
\newcommand{\cycl}{\mathrm{cycl}}

\newcommand{\FT}{\mathrm{FT}}
\newcommand{\Prr}{\mathrm{Pr}}
\newcommand{\restr}{\mathrm{res}}
\newcommand{\C}{\mathbb{C}}
\newcommand{\reg}{\mathrm{reg}}
\newcommand{\gd}{\mathrm{gd}}
\newcommand{\Rbb}{\mathbb{R}}
\newcommand{\Hfrak}{\mathfrak{H}}

\newcommand{\br}{\mathrm{br}}
\newcommand{\cbf}{\mathbf{c}}

\newcommand{\wbf}{\mathbf{w}}
\newcommand{\tbf}{\mathbf{t}}
\newcommand{\sbf}{\mathbf{s}}
\newcommand{\xbf}{\mathbf{x}}

\title[Disjointness of cohomology groups]{On the Broué--Digne--Michel conjectures: disjointness of the cohomology groups}	

\author{Arnaud Eteve}

\begin{document}

\maketitle

\noindent\textbf{Abstract.} The goal of this paper is to start proving the Broué--Digne--Michel conjectures concerning the cohomology of Deligne--Lusztig varieties attached to roots of the full twist. The main result establishes that the $\mathrm{Ext}$-algebras of the cohomology complexes attached to roots of the full twist are concentrated in degree $0$. The main idea is to reduce the problem to a calculation in the Hecke category. The resulting argument is essentially topological and applies to all realizations of the Hecke category. 

\tableofcontents

\section{Introduction}

One of the most famous conjectures in modular representation theory of finite groups is Broué's abelian defect conjecture. 
Let $\Gamma$ be a finite group, $\ell$ a prime dividing the order of $\Gamma$ and $k$ an algebraically closed field of characteristic $\ell$.

\begin{conj}[Abelian defect conjecture, \cite{BroueIsometrieParfaites}]\label{conj:abelian-defect}
Let $b \subset k[\Gamma]$ be a block with abelian defect group $D \subset \Gamma$ and Brauer correspondent $b_{D}$ in $N_{\Gamma}(D)$. Then there exists an equivalence of triangulated categories 
$$\DD^b(b) \cong \DD^b(b_D).$$
\end{conj}

This conjecture was originally proposed to understand the structure of blocks of finite groups of Lie type in non-defining characteristic. Let us fix a connected reductive group $\Gbf$ over $\Fqb$, the algebraic closure of a finite field, and let $\Frob : \Gbf \to \Gbf$ be a Frobenius endomorphism, for example a geometric Frobenius coming from some $\Fq$-structure of $\Gbf$. Let us also fix a Borel pair $\Bbf = \Tbf\Ubf$ which we assume to be stable under $\Frob$ and denote by $\Wbf$ the corresponding Weyl group and $\delta : \Wbf \to \Wbf$ the morphism induced by $\Frob$. In that situation, the modular representation theory of $\Gbf^{\Frob}$ in characteristic $\ell \neq \mathrm{char}(\Fqb)$ is mostly controlled by Deligne--Lusztig theory, originally introduced in \cite{DeligneLusztig}. When $\ell$ is large enough, i.e. not dividing the order of $\Wbf$, it is expected that the functor providing the equivalence of Conjecture \ref{conj:abelian-defect} can be constructed using suitably enhanced Deligne--Lusztig induction functors. We refer to \cite[Chapter VI]{BroueReflectionGroupsLectures} for a discussion.

In this context, Broué and Michel \cite{BroueMichelActionRegularElements}, together with Digne and Michel \cite{DigneMichelParabolic} in the parabolic case, proposed a series of conjectures connecting Deligne--Lusztig theory with braid group actions, cyclotomic Hecke algebras \cite{BroueMalleCyclotomicHeckeAlg} and the theory of generic reductive groups in the sense of \cite{BroueMalle, BroueMalleMichel}.  Precisely, let $\Br_{\Wbf}$ be the braid group of $\Wbf$ and, given $w \in \Wbf$, we denote by $\wbf \in \Br_{\Wbf}$ its canonical lift. There is a distinguished central element 
$$\FT = \wbf_0^2 \in \Br_{\Wbf}$$ 
called the full twist, where $w_0 \in \Wbf$ is the longest element. 

\begin{defi}
An element $c\delta \in \Wbf\delta$ is called a root of the full twist if 
$$(\cbf\delta)^d = \FT\delta^d$$
and $c\delta$ is a good root of the full twist if 
$$\ell((c\delta)^i) = i\ell(c\delta), 1 \leq i \leq \left\lfloor \frac{d}{2}\right\rfloor.$$ 
\end{defi}

\begin{conj}[\protect{\cite[Conjecture 2.2]{DigneMichelEndomorphism}}]\label{conj-intro:good-roots-of-full-twists}
All good roots of the full twists are conjugate under $\Br_{\Wbf}$.
\end{conj}

\begin{thm}[\protect{\cite[Theorem 12.4 (ii)]{Bessis}}]
Assume that $\delta = 1$ (e.g. $\Gbf$ is split), then Conjecture 
\end{thm}

For a fixed value of $d$, the good roots of the full twist are special elements of $\Wbf\delta$. They are all regular in the sense of Springer \cite{SpringerRegularElements}. In particular, the centralizer of a regular element $c\delta$, which we denote by $\Wbf(c\delta)$, is a complex reflection group. 

 Given an element $w \in \Wbf$, we denote by 
$$X(w) := \{g\Bbf \in \Gbf/\Bbf, g^{-1}\Frob(g) \in \Bbf w\Bbf\}$$
the corresponding Deligne--Lusztig variety.

\begin{conj}[\protect{\cite[Conjecture 5.7]{BroueMichelActionRegularElements}}]\label{conj-intro:BroueMichel}
Let $c\delta \in \Wbf\delta$ be a good root of the full twist, then 
\begin{enumerate}
\item The cohomology groups $H^i_c(X(c), \Qlb)$ are disjoint, i.e., for all $i \neq j$, we have 
$$\Hom_{\Gbf^{\Frob}}(H^i_c(X(c), \Qlb), H^j_c(X(c), \Qlb)) = 0.$$
\item There is a canonical isomorphism 
$$\Br_{\Wbf}(\cbf\delta) \cong \Br_{\Wbf(c\delta)}.$$
\item There is a canonical map 
$$\Qlb[\Br_{\Wbf}(\cbf\delta)] \to \End(H^*_c(X(c), \Qlb))$$
which factors through a cyclotomic Hecke algebra of the complex reflection group $\Wbf(c\delta)$. 
\item The groups $H^i_c(X(c), \Qlb)$ decompose in an explicit way in terms of characters of $\Gbf^{\Frob}$ and the corresponding Hecke algebra of $\Wbf(c\delta)$.
\end{enumerate}
\end{conj}

\begin{rques}
\begin{enumerate}
\item These conjectures also provide a collection of examples where we expect a description of the individual cohomology groups of the Deligne--Lusztig varieties $X(w)$. 
\item These conjectures are completely proven when $c\delta$ is a ($\delta$-twisted) Coxeter element by \cite{LusztigCoxeterOrbits}. 
\item Item $(2)$ of the conjecture is a theorem of Bessis \cite{Bessis} when $\delta = 1$.
\item The parabolic analogs of this conjecture are stated in \cite[Section 9]	{DigneMichelParabolic}. 
\item The map $\Qlb[\Br_{\Wbf}^+(\cbf\delta)] \to \End(H^*_c(X(c), \Qlb))$ is already constructed in \cite{BroueMichelActionRegularElements}.
\end{enumerate}
\end{rques}
 
Our main theorems concern the disjointness part of Conjecture \ref{conj-intro:BroueMichel}. We let $\Lambda \in \{\Flb, \Zlb, \Qlb\}$ be a coefficient ring.

\begin{thm}[Theorem \ref{thm:concentration-degree-0}]\label{thm-intro:concentration-degree-0}
There exists $c\delta$ a good root of the full twist such that
$$\End(\RGamma_c(Y(\dot{c}), \Lambda))$$
lies in degree $0$ and in particular over $\Lambda = \Qlb$, 
$$\dim(\End(\RGamma_c(X(c), \Lambda))) = |\Wbf(c\delta)|.$$
\end{thm}

\begin{thm}[Theorem \ref{thm:surjectivity}]\label{thm-intro:surjectivity}
Assume that $\delta = 1$ and let $c \in \Wbf$ be a good root of the full twist, then the map 
$$\Qlb[\Br_{\Wbf}(\cbf)] \to \End(\RGamma_c(X(c), \Qlb))$$
is surjective. 
\end{thm}

\begin{thm}[Theorem \ref{thm:vanishing-non-extendable-elements}]\label{thm-intro:concentration-degree-0-parabolic}
There exists $c\delta \in \Wbf\delta$, a non-extendable $(d,2k)$-periodic element with $2k \leq d$, see Definition \ref{defi:periodic-elements}, such that the endomorphism algebra 
$$\End_{\Gbf^{\Frob}}(\RGamma_c(Y(\dot{c}, \Pbf), \Lambda))$$ 
is in degree $0$. 
\end{thm}

\begin{conj}\label{conj-intro:conj-of-non-extendable-roots}
All non-extendable $(d,2k)$-periodic elements are conjugate under $\Br_{\Wbf}$.
\end{conj}

\begin{rque}
In particular, if Conjecture \ref{conj-intro:conj-of-non-extendable-roots} holds, then Theorem \ref{thm-intro:concentration-degree-0-parabolic} holds for all non-extendable $(d,2k)$-periodic elements.
\end{rque}

Previously known cases of these theorems, with $\Lambda = \Qlb$ and $X(w)$ instead of $Y(\dot{w})$, include the following:  
\begin{enumerate}
\item when $c = w_0$, Theorems \ref{thm-intro:concentration-degree-0} and \ref{thm-intro:surjectivity} are proven by \cite{LusztigRepFiniteChevalleyGroups},
\item when $c$ is a Coxeter element of minimal length (possibly $\delta$-twisted), then Theorems \ref{thm-intro:concentration-degree-0} and \ref{thm-intro:surjectivity} are theorems of \cite{LusztigCoxeterOrbits},
\item when $c = \FT$, Theorems \ref{thm-intro:concentration-degree-0} and \ref{thm-intro:surjectivity} are due to \cite{BonnafeDudasRouquier}, 
\item when $\Gbf$ has split type $A$, Theorems \ref{thm-intro:concentration-degree-0}, \ref{thm-intro:surjectivity} and \ref{thm-intro:concentration-degree-0-parabolic} are due to \cite{DudasCohomologyDLTypeA} and \cite{BonnafeDudasRouquier},
\item for groups of rank less than $2$, Theorems \ref{thm-intro:concentration-degree-0} and \ref{thm-intro:surjectivity} are due to \cite{DigneMichelRouquier}.
\end{enumerate}

\begin{rque}
We plan to address the rest of Conjecture \ref{conj-intro:BroueMichel} together with its application to Conjecture \ref{conj:abelian-defect} in subsequent work. 
\end{rque}

Let us now briefly comment on the proof of Theorem \ref{thm-intro:concentration-degree-0}. The proof of Theorem \ref{thm-intro:concentration-degree-0-parabolic} is a variation on the same set of ideas and the proof of Theorem \ref{thm-intro:surjectivity} relies on previously known calculations, chiefly those in \cite{DigneMichelRouquier}. Let $\Hcal_{\Lambda}$ be the free monodromic Hecke category in the sense of \cite{BezrukavnikovYun, BezrukavnikovRicheSoergelTheory,Gouttard, EteveFreeMono}. We review some basic aspects of it in Section \ref{sec:free-mono-hecke}. For the purpose of this introduction, let us mention two facts: 
\begin{enumerate}
\item there is a collection of invertible objects $(\Delta_w)_w$ called standard objects which satisfy the braid relations and whose inverses $(\nabla_{w^{-1}})_{w}$ are called costandard objects, 
\item there is a natural functor $\Hcal_{\Lambda} \xrightarrow{\ch\pfrak_!} \DD(\Rep_{\Lambda}\Gbf^{\Frob})$, see Section \ref{sec:braid-actions}, and isomorphisms
$$\ch\pfrak_!(\Delta_w) \cong \RGamma_c(Y(\dot{w}), \Lambda)[\ell(w)].$$
\end{enumerate}

The proof of Theorem \ref{thm-intro:concentration-degree-0} then goes as follows: 
\begin{enumerate}
\item we first perform a series of geometric reductions. It is enough to show that for all $x \in \Wbf$, the complex 
\begin{equation}\label{eq-intro:mapping-space}
\Hom_{\Hcal_{\Lambda}}(\Delta_x\Delta_c, \Delta_c\nabla_{\delta(x)})
\end{equation}
is concentrated in degree $0$.
\item We then use some positivity properties in the Hecke algebra, which translate into the following statement in the Hecke category: if the mapping space \eqref{eq-intro:mapping-space} is nonzero, then we must have $xc \leq c\delta(x)$ and $\delta(x)c^{-1} \leq c^{-1}x$, this is Lemma \ref{lem:control-mapping-space} and whenever $xc = c\delta(x)$, then \eqref{eq-intro:mapping-space} is isomorphic to $\RR_{\Lambda}[0]$.
\item Since all good roots of the full twist are conjugate to each other in the braid group, we are free to pick any one of them to prove the statement, we show that we can pick such a good root of the full twist $c\delta$, such that whenever $xc \leq c\delta(x)$ and $\delta(x)c^{-1} \leq c^{-1}x$, then $xc = c\delta(x)$, this is Lemma \ref{lem:good-root-of-ft-selects-centralizer}.
\end{enumerate}

\subsection{Notation and conventions}

\subsubsection{Reductive groups} We let $\Gbf$ be a reductive group over $\Fqb$ which comes equipped with a Frobenius endomorphism 
$$ \Frob : \Gbf \to \Gbf$$ 
and a $\Frob$-stable Borel pair $\Bbf = \Tbf\Ubf$. We denote by $\Wbf$ the corresponding Weyl group, by $\delta : \Wbf \to \Wbf$ the morphism induced by $\Frob$ and by $\widetilde{\Wbf} = \Wbf \rtimes \langle \delta \rangle$. 

\subsubsection{Root datum}
We denote by $(X,X^{\vee}, \Phi, \Phi^{\vee}, \Delta, \Delta^{\vee})$ the based root datum attached to $(\Gbf, \Bbf, \Tbf)$ and by $S \subset \Wbf$ the corresponding set of simple reflections. We have in particular that $\delta(S) = S$. For any $w\in \Wbf$, we denote by $\dot{w} \in \Nbf_{\Gbf}(\Tbf)$ the Tits lift of $w$. For an element $n \in \Nbf_{\Gbf}(\Tbf)$ mapping to $w \in \Wbf$, we shall write $\ell(n)= \ell(w)$.

\subsubsection{Parabolic subgroups}

Given $I \subset S$ a subset of simple reflections, we denote by $\Pbf = \Lbf\UP$ the corresponding standard parabolic and Levi subgroups. We denote by $\Wbf_I = \Wbf_{\Pbf}$ the Weyl group of $(\Lbf, \Tbf)$. Given two standard parabolic subgroups $\Pbf, \Qbf$, we denote by ${^{\Pbf}}\Wbf^{\Qbf}$ the minimal double coset representatives of $\Wbf_{\Pbf}\backslash \Wbf/\Wbf_{\Qbf}$.

\subsubsection{Bruhat ordering}\label{ssec:bruhat-ordering}

We denote by $\leq$ the Bruhat order on $\Wbf$, we recall, see \cite[Section 5.9]{HumphreysReflectionGroups}, that this is the transitive closure of the following relation
$w' \to w$ if there exists a roots $\alpha \in \Phi$ such that $w' = ws_{\alpha}$ and $\ell(w') \leq \ell(w)$. Moreover if $\alpha > 0$, if $w' < w$ and $w = s_{\alpha}w'$, then $w^{',-1}\alpha > 0$. In particular, given $x < y$ two elements of $\Wbf$, there exists a sequence $x_i$, $i = 1, \dots, n$ with $x_1 = x, x_n = y$ and roots $\alpha_i > 0$ such that 
\begin{enumerate}
\item $x_{i+1} = s_{\alpha_i}x_i$, 
\item $x_i^{-1}\alpha_i > 0$.
\end{enumerate}

\subsubsection{Root systems}

We denote by $V = V_{\Rbb} = X \otimes_{\Z} \Rbb$ the real realization of the root datum and set $V_{\C} = V \otimes_{\Rbb} \C$. We denote by $V_{\C, \reg} \subset V_{\C}$ the regular locus, i.e. the complement of reflection hyperplanes attached to reflections in $\Wbf$. We denote by $C^+ \subset V_{\reg}$ the dominant chamber, by $\overline{C}^+$ its closure, and for $I \subset S$, we denote by $\overline{C}^+_I$ the closure of the facet corresponding to $I$. We denote by $\Hfrak$ the set of reflecting hyperplanes in $V$ (or equivalently in $V_{\C}$) and for any convex subset $K \subset V_{\C}$, we denote by $\Hfrak_K \subset \Hfrak$ the set of hyperplanes that contain $K$. An element $k \in K$ is regular if $k \in K \setminus \bigcup_{H \not \in \Hfrak_K} H$. We denote by $K_{\reg}$ the subset of regular elements in $K$. We fix a $\Wbf$-stable scalar product on $V$ which we denote by $\langle - ,- \rangle$.

\subsubsection{Braid groups}

We denote by $\Br_{\Wbf}^+$ and $\Br_{\Wbf}$ the corresponding braid monoid and braid group. We denote by 
$$\Wbf \to \Br_{\Wbf}^+, w \mapsto \wbf$$
the canonical set-theoretic, length-additive section of the map $\Br^+_{\Wbf} \to \Wbf$ which maps $s \in S$ onto the corresponding generator of $\Br_{\Wbf}^+$. The action of $\delta$ on $\Wbf$ lifts to $\Br_{\Wbf}^+$ and $\Br_{\Wbf}$ and we denote by $\widetilde{\Br}^+_{\Wbf} = \Br_{\Wbf}^+ \rtimes \langle \delta \rangle$ and $\widetilde{\Br}_{\Wbf} = \Br_{\Wbf} \rtimes \langle \delta \rangle$. 

We denote by $\FT \in \Br_{\Wbf}^+$ the full twist, i.e. 
$$\FT = \wbf_0^2,$$
where $w_0 \in \Wbf$ is the longest element. For an element $w \in \Wbf$, we denote its centralizer by $\Wbf(w)$ and similarly, for an element $b \in \Br_{\Wbf}$, we denote its centralizer by $\Br_{\Wbf}(b)$.

\subsubsection{Sheaves} We set $p = \operatorname{char}(\Fq)$ and fix a prime $\ell \neq p$ and denote by $\Lambda \in \{\Flb,\Zlb, \Qlb\}$ a coefficient ring. For an algebraic stack $X$ defined over $\Fqb$, we denote by $\DD(X, \Lambda)$ (or $\DD(X)$ if $\Lambda$ is clear from the context) the category of ind-constructible sheaves on $X$ as defined in \cite{GaitsgoryLurie}.

\subsubsection{Categorical traces and categories}
We recall some basic definitions and notation about categorical traces. We refer to \cite{ToyModel, ZhuCategoricalTame, HoyoisScherotzkeSibilla} for some general background. Assume that we have fixed a base ring $\Lambda$. We denote by $\Prr_{\Lambda}$ the category of presentable $\Lambda$-linear stable categories equipped with the Lurie tensor product. For any dualizable category $\Ccal \in \Prr_{\Lambda}$ and any continuous endomorphism $f : \Ccal \to \Ccal$, we denote by 
$$\Tr(f, \Ccal) \in \DD(\Lambda)$$
the trace of the endomorphism $f$ on the object $\Ccal$. For a monoidal category $\Acal \in \Alg(\Prr_{\Lambda})$, we set 
$$\Prr_{\Acal} := \Mod_{\Acal}(\Prr_{\Lambda}).$$
A continuous endomorphism of $\Prr_{\Acal}$ is induced by an $\Acal$-bimodule $\Mcal$ and we denote by 
$$[\Mcal] : \Prr_{\Acal} \to \Prr_{\Acal}, \Ncal \mapsto \Mcal \otimes_{\Acal} \Ncal$$ 
the induced endomorphism. 
We write
$$\Tr([\Mcal], \Prr_{\Acal})$$
for the trace of the morphism $[\Mcal]$ on the presentable $(\infty, 2)$-category $\Prr_{\Acal}$ (which is viewed as a dualizable object in the $(\infty,3)$-category of presentable $(\infty, 2)$-categories, in the sense of Scholze--Stefanich \cite{ScholzeStefanich, Stefanich}). Recall that this is the Hochschild homology 
$$\Tr([\Mcal], \Prr_{\Acal}) = \Acal \otimes_{\Acal \otimes_{\Lambda} \Acal^{\rev}} \Mcal$$
where $\Acal^{\rev}$ is the category $\Acal$ with reversed tensor product, i.e. $x \otimes^{\rev} y = y \otimes x$ so that a right $\Acal$-module is the same as a left $\Acal^{\rev}$-module. 

\subsection{Acknowledgments}
The author thanks Jean-François Dat, Olivier Dudas, Cédric Bonnafé, Jean Michel, Dennis Gaitsgory, Wyatt Reeves, Jens Eberhardt, Alexander Ivanov, Kazuma Ohara, Sean Cotner, and Thibaud van den Hove for discussions during the preparation of this paper. The author thanks Jean Michel for pointing out a couple of mistakes on a previous version of this text. The author was partially supported by the European Research Council (ERC) under the European Union's Horizon 2020 research and innovation programme (grant agreement No. 950326). The author thanks the Max Planck Institute for its support during the preparation of this paper. 

\subsection{AI Use} We used ChatGPT Astra 6 during the preparation of this paper for spellchecking, proofreading, and literature searches. Lemmas \ref{lem:good-root-of-ft-selects-centralizer}, \ref{lem:lenght-intersection-of-unip-groups} and \ref{lem:coset-centralizer-of-periodic-elements} we suggested and proven by AI. These were obtained during a long discussion about the strategy for Theorem \ref{thm-intro:concentration-degree-0}, we had provided a long chain a possible reductions and during the discussion, the AI suggested that if we add the content of Lemma \ref{lem:good-root-of-ft-selects-centralizer} then the statement was complete. This text has been fully written by the author and the author bears full responsibility for it. 

\section{Preliminaries}

\subsection{The free monodromic unipotent Hecke category}\label{sec:free-mono-hecke}

We denote by $\Hcal_{\Lambda}$ the free monodromic Hecke category as defined in \cite{BezrukavnikovYun, BezrukavnikovRicheSoergelTheory,Gouttard, EteveFreeMono}. This is a certain category of twisted equivariant sheaves on the stack $\Ubf \backslash \Gbf /\Ubf$. Recall that, as in \emph{loc. cit.}, we denote by $\RR_{\Lambda}$ (or simply $\RR$ if $\Lambda$ is clear from context), the ring 
$$\RR_{\Lambda} := \Ocal(\widehat{\Tbf}_{\Lambda})^{\wedge}_{1}$$
where $\widehat{\Tbf}_{\Lambda}$ is the torus dual to $\Tbf$ defined over $\Lambda$ and $\Ocal(\widehat{\Tbf}_{\Lambda})^{\wedge}_{1}$ is the completion of $\Ocal(\widehat{\Tbf}_{\Lambda})$ along the unit section of $\widehat{\Tbf}_{\Lambda}$. We also denote by 
$$\Lcal \in \DD(\Tbf, \RR)$$ 
the free monodromic unipotent local system. 

Recall from \cite{EteveJordan} that we denote by $\Ch_{\Lambda}(\Tbf)$ the set of characters $\pi_1^{t}(\Tbf, 1) \to \Lambda^{\times}$ of finite order invertible in $\Lambda$ and that, after choosing a topological generator $\gamma \in \pi_1^t(\Gm, 1) \cong \widehat{\Z}^{(p)}(1)$ we have a one-to-one correspondence 
$$\Ch_{\Lambda}(\Tbf) \cong \widehat{\Tbf}_{\Lambda}(\Lambda)_{f}$$
with the set of $\Lambda$-points of $\widehat{\Tbf}$ of order invertible in $\Lambda$. Given some $s \in \Ch_{\Lambda}(\Tbf)$, there is a well-defined Kummer sheaf $\Lcal_s \in \DD(\Tbf, \Lambda)$ whose monodromy is given by the corresponding character $\pi_1(\Tbf, 1) \to \Lambda^{\times}$.  

For any element $n \in \Nbf_{\Gbf}(\Tbf)$, we denote by 
$$\nu_n : \Nbf_{\Gbf}(\Tbf)_n \cong n\Tbf$$
the trivialization of the connected component of $\Nbf_{\Gbf}(\Tbf)$ containing $n$ induced by $n$. This induces a morphism 
$$\mu_n : \Ubf \backslash \Bbf n \Bbf/\Ubf \to \Nbf_{\Gbf}(\Tbf)_n \cong n\Tbf.$$
Recall that we denote by $\Delta_{n,s}, \nabla_{n,s} \in \Hcal_{\Lambda}$ the standard and costandard objects defined as 
\begin{align*}
\Delta_{n,s} &:= j_{w,!}\mu_n^*(\Lcal \otimes_{\Lambda} \Lcal_s)[\dim(\Tbf) + \ell(w)], \\
\nabla_{n,s} &:= j_{w,*}\mu_n^*(\Lcal \otimes_{\Lambda} \Lcal_s)[\dim(\Tbf) + \ell(w)].
\end{align*}
We denote by 
$$\Delta_{n} := \bigoplus_{s \in \Ch_{\Lambda}(\Tbf)}\Delta_{n,s}, \nabla_{n} := \bigoplus_{s \in \Ch_{\Lambda}(\Tbf)}\nabla_{n,s}.$$
Note that these are infinite sums; in practice, we can always replace all of $\Ch_{\Lambda}(\Tbf)$ by a finite subset, specifically, we can always assume that $\Ch_{\Lambda}(\Tbf)$ is the union of all $\Wbf$-orbits of characters that are stable under $\Frob$. 

\begin{lem}[\protect{\cite[Lemma 4.3.2]{EteveFreeMono}}]
Let $n,n'$ be two elements such that $\ell(n) + \ell(n') = \ell(nn')$, then there is a canonical isomorphism 
$$\Delta_{n} * \Delta_{n'} \xrightarrow{\sim} \Delta_{nn'}.$$
\end{lem}

Deligne's $K(\pi, 1)$ theorem gives a uniquely defined monoidal functor, see also \cite[Section 1.2]{BonnafeDudasRouquier} 
$$\Delta : \Br_{\Wbf} \to \Hcal_{\Lambda}, \wbf \mapsto \Delta_{\dot{w}},$$
where we view $\Br_{\Wbf}$ as a discrete monoidal category with monoidal structure induced by multiplication in $\Br_{\Wbf}$. Alternatively, consider $\DD(\Br_{\Wbf}, \Lambda)$, the category of sheaves on the discrete scheme $\Br_{\Wbf}$ equipped with monoidal structure obtained from convolution coming from the group structure of $\Br_{\Wbf}$. The functor $b \mapsto \Delta_b$ defines a canonical monoidal functor 
$$\Delta : \DD(\Br_{\Wbf}, \Lambda) \to \Hcal_{\Lambda}.$$

\begin{defi}
\begin{enumerate}
\item A $\Delta$-filtration on a perverse sheaf $T \in \Hcal_{\Lambda}^{\heartsuit}$ is a filtration whose graded pieces are isomorphic to direct summands of $\Delta_w$.
\item A $\nabla$-filtration on a perverse sheaf $T \in \Hcal_{\Lambda}^{\heartsuit}$ is a filtration whose graded pieces are isomorphic to direct summands of $\nabla_w$. 
\item An object is tilting if it admits both a $\Delta$-filtration and a $\nabla$-filtration.
\end{enumerate}
\end{defi} 

\begin{thm}[\protect{\cite[Lemma 9.5.1]{Gouttard}, \cite[Lemma 5.9]{BezrukavnikovRicheSoergelTheory}}]
For any pair $(w,s)$ with $w \in \Wbf$ and $s\in\Ch_{\Lambda}(\Tbf)$, there exists a unique (up to isomorphism) indecomposable tilting sheaf $T_{w,s}$ satisfying the following properties: 
\begin{enumerate}
\item the sheaf $T_{w,s}$ is supported on the closure of $\Bbf w\Bbf$,
\item the sheaf $T_{w,s}$ contains $\Delta_{w,s}$ with multiplicity one in one (equivalently any) $\Delta$-filtration.
\end{enumerate}
\end{thm}

\begin{lem}[\protect{\cite[Section 9.8.1]{Gouttard}}]\label{lem:structure-simple-reflection}
Let $t$ be a simple reflection and $s \in \Ch_{\Lambda}(\Tbf)$, one of the following is true 
\begin{enumerate}
\item there are canonical isomorphisms 
$$\Delta_{t,s} \cong T_{t,s} \cong \nabla_{t,s}$$
\item there are fiber sequences 
\begin{align*}
&\Delta_{t,s} \to T_{t,s} \to \Delta_{1,s} \\
&\Delta_{1,s} \to T_{t,s} \to \nabla_{t,s} \\
&\Delta_{t,s} \to \nabla_{t,s} \to \Delta_{1,s}/(1 - e^{-\alpha_t}),
\end{align*}
where $e^{-\alpha_t} \in \RR_{\Lambda}$ is the image of $e^{-\alpha_t} \in \Lambda[X_*]$.
\end{enumerate}
\end{lem}

\begin{lem}\label{lem:product-simple-reflection}
Let $s$ be a simple reflection, then the object $\Delta_s * \Delta_s$ has a filtration whose graded pieces are summands of $\Delta_s, \Delta_s[-1], \Delta_1$.
\end{lem}

\begin{proof}
Since $\nabla_s * \Delta_s = \Delta_1$, this follows from Lemma \ref{lem:structure-simple-reflection}.
\end{proof}

\subsection{Deligne--Lusztig theory and braid actions}\label{sec:braid-actions}

\subsubsection{Deligne--Lusztig theory and horocycle correspondence}
Consider the diagram of algebraic stacks 
\[\begin{tikzcd}
	{\Ubf \backslash \Gbf /\Ubf} && \\
	& {\frac{\Gbf}{\Ad_{\Frob}(\Bbf)}} \\
	{\frac{\Ubf \backslash \Gbf /\Ubf}{\Ad_{\Frob}(\Tbf)}} && {\frac{\Gbf}{\Ad_{\Frob}(\Gbf)}}
	\arrow["\pfrak"', from=1-1, to=3-1]
	\arrow["\alpha", from=2-2, to=3-1]
	\arrow["\beta"', from=2-2, to=3-3]
\end{tikzcd}\]

Recall that we denote by 
\begin{enumerate}
\item $\ch = \beta_!\alpha^* : \DD(\frac{\Ubf \backslash \Gbf /\Ubf}{\Ad_{\Frob}(\Tbf)}) \to \DD(\frac{\Gbf}{\Ad_{\Frob}(\Gbf)})$, 
\item $\hc^! = \alpha_!\beta^*[2\dim(\Ubf)] : \DD(\frac{\Gbf}{\Ad_{\Frob}(\Gbf)}) \to \DD(\frac{\Ubf \backslash \Gbf /\Ubf}{\Ad_{\Frob}(\Tbf)})$, 
\item $\hc^* = \alpha_*\beta^! : \DD(\frac{\Gbf}{\Ad_{\Frob}(\Gbf)}) \to \DD(\frac{\Ubf \backslash \Gbf /\Ubf}{\Ad_{\Frob}(\Tbf)}).$ 
\end{enumerate}
These functors form an adjoint triple 
$$(\hc^! \dashv \ch \dashv \hc^*).$$

Recall as well that we have 
\begin{lem}[\protect{\cite[Lemma 3.2.7]{EteveFreeMono}}]\label{lem:dl-vs-ch}
There is a canonical isomorphism
$$\ch \circ \pfrak_!\Delta_{\dot{w}} \cong \RGamma_c(Y(\dot{w}), \Lambda)[\ell(w)]$$
in $\DD(\frac{\Gbf}{\Ad_{\Frob}(\Gbf)})$, where $Y(\dot{w}) = \{g\Ubf \in \Gbf/\Ubf, g^{-1}\Frob(g) \in \Ubf \dot{w}\Ubf\}$ is the usual Deligne--Lusztig variety of \cite{DeligneLusztig}.
\end{lem}

\begin{rque}\label{rque:series}
If we only care about the unipotent part of the cohomology of $Y(\dot{w})$, we can instead use $\Delta_{\dot{w},1}$. In particular, we have 
$$\ch \circ \pfrak_!\Delta_{\dot{w},1} \cong e_1^{\Lambda}\RGamma_c(Y(\dot{w}), \Lambda)[\ell(w)]$$
where $e_1^{\Lambda} \in \Lambda[T^{w\Frob}]$ is the idempotent corresponding to the principal block. Analogous statements hold for other Lusztig series. We refer to \cite[Section 5.3]{EteveJordan} for a discussion.
\end{rque}

\subsubsection{Parabolic character sheaves}

Let $\Pbf = \Lbf\UP$ be a standard parabolic subgroup. We recall the construction of parabolic character sheaves of \cite{LusztigParabolicCS1} (see also \cite[Section 4]{EteveJordan} for a discussion). The category of parabolic character sheaves is defined to be the category of sheaves on the stack 
$$\frac{\UP \backslash \Gbf/\Frob(\UP)}{\Ad_{\Frob}\Lbf}.$$

Let $v \in {^\Pbf}\Wbf^{\Frob(\Pbf)}$ such that $\Ad(\dot{v})\Frob(\Lbf) = \Lbf$, then the choice of $\dot{v}$ yields an isomorphism 
$$\frac{\UP \backslash \Pbf v\Frob(\Pbf)/\Frob(\UP)}{\Ad_{\Frob}\Lbf} \cong \pt/\Lbf^{\dot{v}\Frob} \ltimes (\UP \cap {^{\dot{v}}}\Frob(\UP)),$$
and in particular, there is an equivalence of categories 
$$\DD(\frac{\UP \backslash \Pbf v\Frob(\Pbf)/\Frob(\UP)}{\Ad_{\Frob}\Lbf}) \cong \DD(\Rep_{\Lambda}\Lbf^{\dot{v}\Frob})$$
induced by pullback along $\pt/\Lbf^{\dot{v}\Frob} \ltimes (\UP \cap {^{\dot{v}}}\Frob(\UP)) \to \pt/\Lbf^{\dot{v}\Frob}$. We shall also denote the inclusion by 
$$j_v^{\Pbf} : \frac{\UP \backslash \Pbf v\Frob(\Pbf)/\Frob(\UP)}{\Ad_{\Frob}\Lbf}  \subset  \frac{\UP \backslash \Gbf/\Frob(\UP)}{\Ad_{\Frob}\Lbf}.$$

\begin{rque}
The stack $\frac{\UP \backslash \Gbf/\Frob(\UP)}{\Ad_{\Frob}\Lbf}$ admits a stratification whose strata are all of the form $\pt/(\Lbf')^{\dot{v}\Frob}\ltimes \Vbf$ where $\Lbf'\subset \Lbf$ are certain standard Levi subgroups and $\Vbf$ are connected unipotent groups. We refer to \cite[Section 4.2]{EteveJordan} for a discussion. 
\end{rque}

Given $\Qbf = \Mbf\UQ$ a second standard parabolic with $\Pbf \subset \Qbf$, we have a horocycle correspondence 
$$\frac{\UP \backslash \Gbf/\Frob(\UP)}{\Ad_{\Frob}\Lbf} \xleftarrow{\alpha_{\Pbf,\Qbf}} \frac{\UQ \backslash \Gbf/\Frob(\UQ)}{\Ad_{\Frob}(\Mbf\cap \Pbf)} \xrightarrow{\beta_{\Pbf, \Qbf}} \frac{\UQ \backslash \Gbf/\Frob(\UQ)}{\Ad_{\Frob}\Mbf}$$
and an adjoint triple 
$$(\hc_{\Pbf,\Qbf}^! \dashv \ch_{\Pbf, \Qbf} \dashv \hc^*_{\Pbf, \Qbf})$$
defined as 
\begin{enumerate}
\item $\hc_{\Pbf,\Qbf}^! = \alpha_{\Pbf, \Qbf,!}\beta_{\Pbf,\Qbf}^*[2\dim(\UP/\UQ)]$
\item $\ch_{\Pbf,\Qbf} = \beta_{\Pbf, \Qbf,!}\alpha_{\Pbf,\Qbf}^*$,
\item $\hc_{\Pbf,\Qbf}^* = \alpha_{\Pbf, \Qbf,*}\beta_{\Pbf,\Qbf}^!$.
\end{enumerate}  
We shall denote by 
$$\hc_{\Pbf}^? = \hc_{\Pbf,\Gbf}^?, \ch_{\Pbf} = \ch_{\Pbf,\Gbf}$$
and by 
$$\hc^{\Pbf,?} = \hc_{\Bbf,\Pbf}^?, \ch^{\Pbf} = \ch_{\Bbf,\Pbf},$$
where $?\in\{*,!\}$.

\begin{lem}[\protect{\cite[Lemma 3.2.7]{EteveFreeMono}}]\label{lem:comparison-with-dl-parabolic}
There is a natural isomorphism of functors $\DD(\Rep_{\Lambda}\Lbf^{\dot{v}\Frob}) \to \DD(\Rep_{\Lambda}\Gbf^{\Frob})$ 
$$\ch_{\Pbf}j_{v,!}^{\Pbf} \cong \RGamma_c(Y(\dot{v}, \Pbf), \Lambda) \otimes_{\Lambda[\Lbf^{\dot{v}\Frob}]} - $$
where $Y(\dot{v}, \Pbf)$ is the usual parabolic Deligne--Lusztig variety defined as 
$$Y(\dot{v}, \Pbf) := \{g\UP \in \Gbf/\UP, g^{-1}\Frob(g) \in \UP \dot{v}\Frob(\UP)\}.$$
\end{lem}

\begin{lem}[\protect{\cite[Lemma 4.1.4]{EteveJordan}}]
For three parabolic subgroups $\Pbf \subset \Qbf \subset \Rbf$, there are natural transitivity isomorphisms 
\begin{enumerate}
\item $\hc_{\Pbf,\Qbf}^!\hc_{\Qbf,\Rbf}^! \cong \hc_{\Pbf,\Rbf}^!$,
\item $\ch_{\Qbf,\Rbf}\ch_{\Pbf,\Qbf} \cong \ch_{\Pbf,\Rbf}$,
\item $\hc_{\Pbf,\Qbf}^*\hc_{\Qbf,\Rbf}^* \cong \hc_{\Pbf,\Rbf}^*$.
\end{enumerate}
\end{lem}

\subsubsection{Categorical trace of Frobenius}

Recall that we have 
\begin{thm}[\protect{\cite[Theorem 1.3.1]{EteveFreeMono}}, \cite{ZhuCategoricalTame},\cite{GaitsgoryRozemblyumVarshavsky}]
There is a canonical isomorphism 
$$\Tr(\restr_{\Frob^*}, \Prr_{\Hcal_{\Lambda}}) \cong \DD(\frac{\Gbf}{\Ad_{\Frob}(\Gbf)})$$ 
and the canonical trace map 
$$\tr_{\Frob} : \Hcal_{\Lambda} \to \DD(\frac{\Gbf}{\Ad_{\Frob}(\Gbf)})$$ 
is given by the functor $\ch \circ \pfrak_!.$
\end{thm}

\begin{rque}\label{rque:canonical-trace-functor}
The functor $\tr_{\Frob}$ being the canonical trace map means that, in particular, it comes equipped with a collection of functorial isomorphisms
$$\tr_{\Frob}(X * Y) \cong \tr_{\Frob}(\Frob^*(Y) * X).$$
Note that as $\Frob^*$ is invertible, we have $\tr_{\Frob}(X * Y) \cong \tr_{\Frob}(Y * \Frob_*(X))$ functorially in $X,Y$. 
\end{rque} 

Consider the category $\DD(\Br_{\Wbf}, \Lambda)$ as a monoidal category equipped with an endofunctor 
$$\Frob^* : \DD(\Br_{\Wbf}, \Lambda) \to \DD(\Br_{\Wbf}, \Lambda)$$ 
induced by $\delta : \Br_{\Wbf} \to \Br_{\Wbf}$. 
\begin{prop}
There is a canonical equivalence 
$$\Tr(\restr_{\Frob^*}, \Prr_{\DD(\Br_{\Wbf}, \Lambda)}) \cong \DD(\frac{\Br_{\Wbf}}{\Ad_{\delta}\Br_{\Wbf}}, \Lambda).$$
\end{prop}

\begin{proof}
We can directly compute 
\begin{align*}
\Tr(\restr_{\Frob^*}, \Prr_{\DD(\Br_{\Wbf}, \Lambda)}) &= \DD(\Br_{\Wbf}) \otimes_{\DD(\Br_{\Wbf}) \otimes \DD(\Br_{\Wbf}^{\opp})} \DD(\Br_{\Wbf})_{\delta} \\
&\cong \DD(\Delta(\Br_{\Wbf}) \backslash (\Br_{\Wbf} \times \Br_{\Wbf})/\Delta_{\delta}\Br_{\Wbf}) \\
&\cong \DD(\frac{\Br_{\Wbf}}{\Ad_{\delta}\Br_{\Wbf}}),
\end{align*}
where the equivalences are as follows:
\begin{enumerate}
\item the first one is that $\Tr(\restr_{\Frob^*}, \Prr_{\DD(\Br_{\Wbf}, \Lambda)})$ is computed via Hochschild homology, 
\item the second is via categorical Künneth formulas which hold because $\Br_{\Wbf}$ is a discrete scheme, 
\item the third one is simply rewriting the quotient.
\end{enumerate}
\end{proof}

The functor 
$$\Delta : \DD(\Br_{\Wbf}, \Lambda) \to \Hcal_{\Lambda}$$ 
is visibly compatible with $\Frob$ acting on both sides, and both categories are rigid. For $\Hcal_{\Lambda}$ this is \cite[Theorem 4.4.2]{EteveFreeMono} and for $\DD(\Br_{\Wbf}, \Lambda)$ this is immediate. It follows from the functoriality of categorical traces, see \cite[Section 3.6]{ToyModel}, that there is a well-defined map 
$$\Tr(\Delta) : \Tr(\restr_{\Frob^*}, \Prr_{\DD(\Br_{\Wbf}, \Lambda)}) \to \Tr(\restr_{\Frob^*}, \Prr_{\Hcal_{\Lambda}}).$$

Let $b,b' \in \Br_{\Wbf}$ be two elements such that $bb' = b'\delta(b)$, then there is a well-defined endomorphism, see \cite[Section 4.1]{BonnafeDudasRouquier}, $D_b^{\tr} : \tr(\Delta_{b'}) \to \tr(\Delta_{b'})$ given by the following composition
\begin{align*}
\tr_{\Frob}(\Delta_{b'}) &\cong \tr_{\Frob}(\Delta_{b}\Delta_{b^{-1}}\Delta_{b'}) \\
&\cong \tr_{\Frob}(\Delta_{b^{-1}}\Delta_{b'}\Frob_*(\Delta_b)) \\
&\cong \tr_{\Frob}(\Delta_{b^{-1}}\Delta_{b'}\Delta_{\delta(b)}) \\
&\cong \tr_{\Frob}(\Delta_{b^{-1}}\Delta_{b}\Delta_{b'}) \\
&\cong \tr_{\Frob}(\Delta_{b'}), 
\end{align*}
where the maps are defined as follows:
\begin{enumerate}
\item the first and last maps are induced by the isomorphism $\Delta_{b}\Delta_{b^{-1}} \cong \Delta_1 \cong \Delta_{b^{-1}}\Delta_{b}$, 
\item the second map is the $\Frob_*$-cyclicity, see Remark \ref{rque:canonical-trace-functor}, of the functor $\tr_{\Frob}$, 
\item the third map is the natural isomorphism $\Frob_*(\Delta_b) \cong \Delta_{\delta(b)}$, 
\item the fourth map is induced by the commutation relation $bb' = b'\delta(b)$. 
\end{enumerate}

Consider the functor $\Tr(\Delta) : \DD(\frac{\Br_{\Wbf}}{\Ad_{\delta}\Br_{\Wbf}}) \to \DD(\frac{\Gbf}{\Ad_{\Frob}(\Gbf)})$, and let $b,b' \in \Br_{\Wbf}$ be two elements such that $bb' = b'\delta(b)$, then evaluation at $b'$ defines a functor
$$\DD(\{b'\}/C_{\Br_{\Wbf}}(b'\delta)) \to \DD(\frac{\Gbf}{\Ad_{\Frob}(\Gbf)})$$
and in particular the object $\tr(\Delta_{b'})$ carries an action of $C_{\Br_{\Wbf}}(b'\delta)$, for $b$ as above, let us denote by $D_b^{\mathrm{act}} : \tr(\Delta_{b'}) \to \tr(\Delta_{b'})$ the action of $b$. 

Finally for $b,b' \in \Br_{\Wbf}^+$ such that $bb' = b'\delta(b)$, \cite{BroueMichelActionRegularElements} constructed an endomorphism 
$$D_b^{\geo} : \tr(\Delta_{bb'}) \to \tr(\Delta_{b'\delta(b)})$$
coming from the first partial Frobenius on the cohomology of Deligne--Lusztig varieties. It is given as follows. Write $b = s_1\dots s_n$ and $b' = c_1\dots c_m$ as products of simple reflections. Recall that the generalized Deligne--Lusztig varieties are defined as 
$$Y(b) := \{g_1\Ubf, \dots, g_n\Ubf \in \Gbf/\Ubf, g_i^{-1}g_{i+1} \in \Ubf \dot{s}_i \Ubf\ (1\leq i<n), g_n^{-1}\Frob(g_1) \in \Ubf\dot{s}_n\Ubf \}.$$
There is a well-defined partial Frobenius 
\begin{align*}
\Frob_b : Y(bb') &\to Y(b'\delta(b)) \\
(g_1\Ubf, \dots, g_{n+m}\Ubf) &\mapsto (g_{n+1}\Ubf, \dots, g_{n+m}\Ubf, \Frob(g_1)\Ubf, \dots, \Frob(g_n)\Ubf).
\end{align*}

\begin{lem}
\begin{enumerate}
\item Let $b,b' \in \Br_{\Wbf}$ such that $bb' = b'\delta(b)$, then we have 
$$D_b^{\tr} = D^{\mathrm{act}}_b.$$
\item Let $b,b' \in \Br_{\Wbf}^+$ such that $bb' = b'\delta(b)$ and set $c = bb'$, then we have 
$$D_b^{\tr} = D_b^{\geo} : \tr(\Delta_{bb'}) \to \tr(\Delta_{b'\delta(b)})$$
where $D_b^{\tr}$ is defined through the equation $b^{-1}c\delta(b) = b'\delta(b)$.
\end{enumerate}
\end{lem} 

\begin{proof}
The comparison $D^{\mathrm{act}}_b = D^{\tr}_b$ is immediate by the functorial properties of categorical traces. Let us now prove the second comparison. Let us introduce the following notation 
\begin{enumerate}
\item $\br : \Delta_{bb'} \xrightarrow{\sim} \Delta_b \Delta_{b'}$ be the braid isomorphism coming from the functor $\Br_{\Wbf} \to \Hcal$,
\item $\unit_l : \Delta_1 \xrightarrow{\sim} \Delta_b \Delta_{b^{-1}}, \unit_r : \Delta_1 \xrightarrow{\sim} \Delta_{b^{-1}} \Delta_b$, 
\item $\counit : \Delta_{b^{-1}}\Delta_b \xrightarrow{\sim} \Delta_1$.
\end{enumerate}
Also observe that there is a natural commutative diagram in $\Hcal$ 
\[\begin{tikzcd}
	{\Delta_{bb'}} & {\Delta_b\Delta_{b^{-1}}\Delta_{bb'}} \\
	{\Delta_{b}\Delta_{b'}} & {\Delta_b\Delta_{b^{-1}}\Delta_{b}\Delta_{b'}}
	\arrow["{\unit_l \cdot bb'}", from=1-1, to=1-2]
	\arrow["\br"', from=1-1, to=2-1]
	\arrow["\br", from=1-2, to=2-2]
	\arrow["{b\cdot \unit_r \cdot b'}"', from=2-1, to=2-2]
\end{tikzcd}\]
where the horizontal maps are induced by multiplying the map $\unit_l$ by $\Delta_{bb'}$ on top and by multiplying the map $\unit_r$ by $\Delta_b (-) \Delta_{b'}$ on the bottom. 
We then have a commutative diagram 
\[\begin{tikzcd}
	{\tr(\Delta_{bb'})} & {\tr(\Delta_b\Delta_{b^{-1}}\Delta_{bb'})} & {\tr(\Delta_{b^{-1}}\Delta_{bb'}\Delta_{\delta(b)})} \\
	{\tr(\Delta_b\Delta_{b'})} & {\tr(\Delta_b\Delta_{b^{-1}}\Delta_b\Delta_{b'})} & {\tr(\Delta_{b^{-1}}\Delta_b\Delta_{b'}\Delta_{\delta(b)})} \\
	{\tr(\Delta_{b'}\Delta_{\delta(b)})} & {\tr(\Delta_{b^{-1}}\Delta_b\Delta_{b'}\Delta_{\delta(b)})} & {\tr(\Delta_{b'}\Delta_{\delta(b)})}
	\arrow["{\unit_l \cdot bb'}", from=1-1, to=1-2]
	\arrow["\br"', from=1-1, to=2-1]
	\arrow["\cycl", from=1-2, to=1-3]
	\arrow["\br"', from=1-2, to=2-2]
	\arrow["\br", from=1-3, to=2-3]
	\arrow["{b\cdot\unit_r \cdot b'}", from=2-1, to=2-2]
	\arrow["\cycl"', from=2-1, to=3-1]
	\arrow["\cycl", from=2-2, to=2-3]
	\arrow["\cycl"', from=2-2, to=3-2]
	\arrow["\counit", from=2-3, to=3-3]
	\arrow["{b\cdot\unit_r \cdot b'}"', from=3-1, to=3-2]
	\arrow[equals, nfold, from=3-2, to=2-3]
	\arrow["\counit"', from=3-2, to=3-3]
\end{tikzcd}\]
where the maps labeled $\cycl$ are obtained by applying the cyclicity operations on traces and using $\Frob_*(\Delta_b) \cong \Delta_{\delta(b)}$.
It follows that the composition along the top and right of the diagram is the operator $D_b^{\tr}$, the composition along the left side of the diagram is the operator $D_b^{\geo}$ and the composition of the maps along the bottom row is the identity. 
\end{proof}

\subsection{Intertwining operators}

Recall from \cite[Section 4.3]{EteveJordan} that we define the intertwining operators as follows 
$$\Int_w^!: \DD(\frac{\Ubf \backslash \Gbf /\Ubf}{\Ad_{\Frob}(\Tbf)}) \to \DD(\frac{\Ubf \backslash \Gbf /\Ubf}{\Ad_{\Frob}(\Tbf)})$$
where $\Int_w^! = \beta_{w,!}\Ad_{\Frob}(\dot{w})_!\alpha_{w}^*$, where the maps are defined as follows 
\[\begin{tikzcd}
	{\frac{\Gbf}{\Ad_{\Frob}(\Bbf \cap {^{w^{-1}}\Bbf})}} &&& {\frac{\Gbf}{\Ad_{\Frob}(\Bbf \cap {^w\Bbf})}} \\
	{\frac{\Ubf \backslash \Gbf /\Ubf}{\Ad_{\Frob}(\Tbf)}} &&& {\frac{\Ubf \backslash \Gbf /\Ubf}{\Ad_{\Frob}(\Tbf)}}
	\arrow["{\Ad_{\Frob}(\dot{w})}"', from=1-1, to=1-4]
	\arrow["{\alpha_w}", from=1-1, to=2-1]
	\arrow["{\beta_w}"', from=1-4, to=2-4]
\end{tikzcd}\]

\begin{lem}[\protect{\cite[Lemma 4.3.2]{EteveJordan}}]\label{lem:mackayfiltration}
The functor $\hc^! \ch$ has a $\Wbf$-indexed filtration whose graded pieces are $\Int_w^!$.
\end{lem}

\begin{lem}\label{lem:intertwinning-operators}
There is a natural isomorphism of functors $\Hcal_{\Lambda} \to \DD(\frac{\Ubf \backslash \Gbf /\Ubf}{\Ad_{\Frob}(\Tbf)})$ 
$$\Int_w^! \pfrak_! \cong \pfrak_!(\Delta_{w^{-1}} * - * \Delta_{\delta(w^{-1})})[-2\dim(\Ubf)].$$ 
\end{lem}

\begin{proof}
Let us first consider the following diagram 
\[\begin{tikzcd}
	& {\frac{\Gbf}{\Ad_{\Frob}(\Ubf \cap {^{w^{-1}}\Ubf})}} & {\frac{\Gbf}{\Ad_{\Frob}(\Ubf \cap {^{w}\Ubf})}} & \\
	{\Ubf \backslash\Gbf/\Ubf} & {\frac{\Gbf}{\Ad_{\Frob}(\Bbf \cap {^{w^{-1}}\Bbf})}} & {\frac{\Gbf}{\Ad_{\Frob}(\Bbf \cap {^{w}\Bbf})}} & {\Ubf \backslash\Gbf/\Ubf} \\
	{\frac{\Ubf \backslash \Gbf /\Ubf}{\Ad_{\Frob}(\Tbf)}} &&& {\frac{\Ubf \backslash \Gbf /\Ubf}{\Ad_{\Frob}(\Tbf)}}
	\arrow["{\Ad_{\Frob}(\dot{w})}", from=1-2, to=1-3]
	\arrow["{\tilde{\alpha}_w}"', from=1-2, to=2-1]
	\arrow[from=1-2, to=2-2]
	\arrow[from=1-3, to=2-3]
	\arrow["{\tilde{\beta}_w}", from=1-3, to=2-4]
	\arrow["\pfrak"', from=2-1, to=3-1]
	\arrow["{\Ad_{\Frob}(\dot{w})}", from=2-2, to=2-3]
	\arrow["{\alpha_w}", from=2-2, to=3-1]
	\arrow["{\beta_w}"', from=2-3, to=3-4]
	\arrow["\pfrak", from=2-4, to=3-4]
\end{tikzcd}\]
as all three squares are Cartesian, it is enough to prove that 
$$\tilde{\beta}_{w,!}\Ad_{\Frob}(\dot{w})_!\tilde{\alpha}_w^* \cong \Delta_{w^{-1}} * - * \Delta_{\delta(w^{-1})}[-2\dim(\Ubf)].$$

Let us rewrite $\Ubf\backslash \Gbf /\Ubf$ as $\Delta(\Gbf) \backslash (\Gbf/\Ubf \times \Gbf/\Ubf)$ and $\frac{\Gbf}{\Ad_{\Frob}(\Ubf \cap {^{w^{-1}}\Ubf})} = \Delta(\Gbf) \backslash (\Gbf \times \Gbf) /\Delta_{\Frob}(\Ubf \cap {^{w^{-1}}\Ubf})$, we can rewrite the top correspondence as 
\[\begin{tikzcd}
	{\Delta(\Gbf) \backslash (\Gbf \times \Gbf) /\Delta_{\Frob}(\Ubf \cap {^{w^{-1}}\Ubf})} & {\Delta(\Gbf) \backslash (\Gbf \times \Gbf) /\Delta_{\Frob}(\Ubf \cap {^{w}\Ubf})} \\
	{\Delta(\Gbf) \backslash (\Gbf/\Ubf \times \Gbf/\Ubf)} & {\Delta(\Gbf) \backslash (\Gbf/\Ubf \times \Gbf/\Ubf)}
	\arrow["{r_{\dot{w}^{-1}, \delta(\dot{w})^{-1}}}", from=1-1, to=1-2]
	\arrow["{\tilde{\alpha}_w}"{description}, from=1-1, to=2-1]
	\arrow["{\tilde{\beta}_w}"{description}, from=1-2, to=2-2]
\end{tikzcd}\]
where $r_{\dot{w}^{-1}, \delta(\dot{w})^{-1}}$ is given by right multiplication by $\dot{w}^{-1}$ and $\delta(\dot{w})^{-1}$ on each component. 
We can further factor this correspondence as follows 
\[\resizebox{0.98\linewidth}{!}{\begin{tikzcd}[column sep=0.5em, ampersand replacement=\&]
	{\Delta(\Gbf) \backslash (\Gbf \times \Gbf) /\Delta_{\Frob}(\Ubf \cap {^{w^{-1}}\Ubf})} \& {\Delta(\Gbf) \backslash (\Gbf \times \Gbf) /\Delta_{\Frob}(\Ubf \cap {^{w}\Ubf})} \\
	{\Delta(\Gbf) \backslash (\Gbf/(\Ubf \cap {^{w^{-1}}\Ubf}) \times \Gbf/(\Ubf \cap {^{w^{-1}}\Ubf}))} \& {\Delta(\Gbf) \backslash (\Gbf/(\Ubf \cap {^{w}\Ubf}) \times \Gbf/(\Ubf \cap {^{w}\Ubf}))} \\
	{\Delta(\Gbf) \backslash (\Gbf/\Ubf \times \Gbf/\Ubf)} \& {\Delta(\Gbf) \backslash (\Gbf/\Ubf \times \Gbf/\Ubf)}
	\arrow["{r_{\dot{w}^{-1}, \delta(\dot{w})^{-1}}}", from=1-1, to=1-2]
	\arrow["{\tilde{\alpha}_w^{(1)}}"', from=1-1, to=2-1]
	\arrow["{\tilde{\beta}_w^{(1)}}", from=1-2, to=2-2]
	\arrow["{r_{\dot{w}^{-1}, \delta(\dot{w})^{-1}}}"', from=2-1, to=2-2]
	\arrow["{\tilde{\alpha}_w^{(2)}}"', from=2-1, to=3-1]
	\arrow["{\tilde{\beta}_w^{(2)}}"', from=2-2, to=3-2]
\end{tikzcd}}\]
where $\tilde{\alpha}_w^{(2)}\tilde{\alpha}_w^{(1)} = \tilde{\alpha}_w$ and similarly for $\tilde{\beta}$. As the middle square is Cartesian and the maps $\tilde{\alpha}_w^{(1)}$ and $\tilde{\beta}_w^{(1)}$ are unipotent of dimension $\dim(\Ubf \cap {^{w}\Ubf})$, we have 
$$\tilde{\beta}_{w,!}\Ad_{\Frob}(\dot{w})_!\tilde{\alpha}_w^* \cong \tilde{\beta}_{w,!}^{(2)}r_{\dot{w}^{-1}, \delta(\dot{w})^{-1},!}\tilde{\alpha}_w^{(2),*}[-2\dim(\Ubf \cap {^{w}\Ubf})].$$ 

We are therefore reduced to proving that for a single copy of $\Gbf$, the correspondence 
\[\begin{tikzcd}
	{\Gbf/ (\Ubf \cap {^{w^{-1}}\Ubf})} & {\Gbf/(\Ubf \cap {^{w}\Ubf})} \\
	{\Gbf/\Ubf} & { \Gbf/\Ubf}
	\arrow["{r_{\dot{w}^{-1}}}", from=1-1, to=1-2]
	\arrow["{\tilde{\gamma}_w}"{description}, from=1-1, to=2-1]
	\arrow["{\tilde{\delta}_w}"{description}, from=1-2, to=2-2]
\end{tikzcd}\]
induces an isomorphism 
$$\tilde{\delta}_{w,!}r_{\dot{w}^{-1},!}\tilde{\gamma}_w^* \cong - * \Delta_{w^{-1}}[-\ell(w)].$$
Consider the diagram 
\[\begin{tikzcd}
	& {\Gbf \times^{\Ubf} \Bbf w^{-1}\Bbf/\Ubf} & {\Gbf \times^{\Bbf} \Bbf w^{-1}\Bbf/\Ubf} & {\Gbf/\Ubf } \\
	{\Gbf/\Ubf} & {\Gbf \times^{\Ubf} \Ubf \dot{w}^{-1}\Ubf/\Ubf} \\
	& {\Gbf / (\Ubf \cap {^{w^{-1}}}\Ubf)}
	\arrow["{m'}", from=1-2, to=1-3]
	\arrow["{\pr_1}"{description}, from=1-2, to=2-1]
	\arrow["m", from=1-3, to=1-4]
	\arrow["\iota"{description}, from=2-2, to=1-2]
	\arrow["{\iota'}"{description}, from=2-2, to=1-3]
	\arrow["{\pr_1}"{description}, from=2-2, to=2-1]
	\arrow["{\tilde{\delta}_w r_{\dot{w}^{-1}}}"', from=3-2, to=1-4]
	\arrow["{\tilde{\gamma}_w}"{description}, from=3-2, to=2-1]
	\arrow["\cong"{description}, from=3-2, to=2-2]
\end{tikzcd}\]
where $\pr_1$ is induced by the first projection, $m$ and $m'$ are induced by the multiplication maps and $\iota, \iota'$ are the natural inclusions.
Finally given some monodromic object $A \in \DD(\Gbf/\Ubf)$, there is a natural isomorphism 
$$m'_!(A \boxtimes_{\Lambda} \oplus_{s \in \Ch(\Tbf)} (\Lcal \otimes_{\Lambda} \Lcal_s)[\dim(\Tbf) + \ell(w)]) \cong  \iota^* A[-\dim(\Tbf) + \ell(w)].$$
The lemma then follows from the commutativity of the diagram. 
\end{proof}

\begin{lem}[\protect{\cite[Lemma 2.8.3]{EteveFreeMono}}]\label{lem:reduction-from-h-to-o}
There is a natural isomorphism, for $A,B \in \Hcal_{\Lambda}$, 
$$\Hom_{\frac{\Ubf \backslash \Gbf /\Ubf}{\Ad_{\Frob}(\Tbf)}}(\pfrak_!A, \pfrak_!B) \cong \Hom_{\Hcal_{\Lambda}}(A, B) \otimes_{\RR \otimes \RR} \RR$$ 
where $\RR$ is viewed as an $\RR$-bimodule via the graph of $\Frob^* : \RR \to \RR$.
\end{lem}

\subsection{Intertwining operators in the parabolic case}

In the parabolic case, we introduce the following intertwining operators as in \cite[Section 4.3]{EteveJordan}. Define the following maps 
\[\begin{tikzcd}
	{\frac{\Gbf}{\Ad_{\Frob}(\Pbf \cap {^{w^{-1}}}\Pbf)}} & {\frac{\Gbf}{\Ad_{\Frob}(\Pbf \cap {^{w}}\Pbf)}} \\
	{\frac{\UP \backslash \Gbf/\Frob(\UP)}{\Ad_{\Frob}(\Lbf)}} & {\frac{\UP \backslash \Gbf/\Frob(\UP)}{\Ad_{\Frob}(\Lbf)}}
	\arrow["{\Ad(\dot{w})}", from=1-1, to=1-2]
	\arrow["{\alpha_w^{\Pbf}}"', from=1-1, to=2-1]
	\arrow["{\beta_w^{\Pbf}}", from=1-2, to=2-2]
\end{tikzcd}\]
the intertwining operator
$$\Int_{\Pbf, \Pbf,w}^! : \DD(\frac{\UP \backslash \Gbf/\Frob(\UP)}{\Ad_{\Frob}(\Lbf)}) \to \DD(\frac{\UP \backslash \Gbf/\Frob(\UP)}{\Ad_{\Frob}(\Lbf)})$$
is defined as $\Int_{\Pbf, \Pbf,w}^! := \beta_{w,!}^{\Pbf}\Ad(\dot{w})_!\alpha_{w}^{\Pbf, !}[2\dim(\UP)]$.

\begin{lem}[\protect{\cite[Lemma 4.3.2]{EteveJordan}}]\label{lem:mackay-filtration-parabolic-case}
The functor $\hc^!_{\Pbf}\ch_{\Pbf}$ has a natural filtration whose graded pieces are $\Int_{\Pbf, \Pbf, w}^!$.
\end{lem}

\begin{lem}\label{lem:comparison-mackay}
The functor $\hc^{\Pbf,!}\Int_{\Pbf, \Pbf,w}^!\ch^{\Pbf}$ has a natural filtration whose graded pieces are given by $\Int_x^!$ where $x \in \Wbf_{\Pbf}w\Wbf_{\Pbf}$. 
\end{lem}

\begin{proof}
This lemma is a refinement of the Mackay filtration of \cite[Lemma 3.2.7]{EteveFreeMono}. Let us introduce the following diagram 
\[\begin{tikzcd}
	{\Xcal_{\Bbf, \Bbf}} & {\Xcal_{\Pbf, \Bbf}} & {\frac{\Gbf}{\Ad_{\Frob}\Bbf}} & {\frac{\UP \backslash \Gbf/\Frob(\UP)}{\Ad_{\Frob}\Lbf}} & {\frac{\Ubf\backslash \Gbf/\Ubf}{\Ad_{\Frob}(\Tbf)}} \\
	{\Xcal_{\Bbf, \Pbf}} & {\Xcal_{\Pbf, \Pbf}} & {\frac{\Gbf}{\Ad_{\Frob}\Pbf}} & {\frac{\UP \backslash \Gbf/\Frob(\UP)}{\Ad_{\Frob}\Lbf}} \\
	{\frac{\Gbf}{\Ad_{\Frob}\Bbf}} & {\frac{\Gbf}{\Ad_{\Frob}\Pbf}} & {\frac{\Gbf}{\Ad_{\Frob}\Gbf}} \\
	{\frac{\UP \backslash \Gbf/\Frob(\UP)}{\Ad_{\Frob}\Lbf}} & {\frac{\UP \backslash \Gbf/\Frob(\UP)}{\Ad_{\Frob}\Lbf}} \\
	{\frac{\Ubf\backslash \Gbf/\Ubf}{\Ad_{\Frob}(\Tbf)}}
	\arrow[from=1-1, to=1-2]
	\arrow[from=1-1, to=2-1]
	\arrow[from=1-2, to=1-3]
	\arrow[from=1-2, to=2-2]
	\arrow[from=1-3, to=1-4]
	\arrow[from=1-3, to=2-3]
	\arrow[from=1-4, to=1-5]
	\arrow[from=1-4, to=2-4]
	\arrow[from=2-1, to=2-2]
	\arrow[from=2-1, to=3-1]
	\arrow[from=2-2, to=2-3]
	\arrow[from=2-2, to=3-2]
	\arrow[from=2-3, to=2-4]
	\arrow[from=2-3, to=3-3]
	\arrow[from=3-1, to=3-2]
	\arrow[from=3-1, to=4-1]
	\arrow[from=3-2, to=3-3]
	\arrow[from=3-2, to=4-2]
	\arrow[from=4-1, to=4-2]
	\arrow[from=4-1, to=5-1]
\end{tikzcd}\]
Lemmas \ref{lem:mackay-filtration-parabolic-case} and \ref{lem:mackayfiltration} are induced by stratifications of $\Xcal_{\Pbf,\Pbf}$ and $\Xcal_{\Bbf, \Bbf}$; specifically, as in \emph{loc. cit.}, there is a natural cube with Cartesian faces 
\[\begin{tikzcd}
	{\Xcal_{\Bbf, \Bbf}} &&& {\Xcal_{\Pbf, \Bbf}} \\
	& {\Bbf\backslash\Gbf/\Bbf} & {\Pbf\backslash\Gbf/\Bbf} \\
	& {\Bbf\backslash\Gbf/\Pbf} & {\Pbf\backslash\Gbf/\Pbf} \\
	{\Xcal_{\Bbf, \Pbf}} &&& {\Xcal_{\Pbf, \Pbf}}
	\arrow[from=1-1, to=1-4]
	\arrow[from=1-1, to=2-2]
	\arrow[from=1-1, to=4-1]
	\arrow[from=1-4, to=2-3]
	\arrow[from=1-4, to=4-4]
	\arrow[from=2-2, to=2-3]
	\arrow[from=2-2, to=3-2]
	\arrow[from=2-3, to=3-3]
	\arrow[from=3-2, to=3-3]
	\arrow[from=4-1, to=3-2]
	\arrow[from=4-1, to=4-4]
	\arrow[from=4-4, to=3-3]
\end{tikzcd}\]
and the inverse image of the stratum $\Pbf \backslash \Pbf w \Pbf /\Pbf$ in $\Xcal_{\Pbf, \Pbf}$ gives rise to the correspondence defining $\Int_{\Pbf,\Pbf, w}^!$ (and similarly with $\Xcal_{\Bbf, \Bbf}$). Taking the inverse image of $\Pbf \backslash \Pbf w \Pbf /\Pbf$ in $\Xcal_{\Bbf, \Bbf}$ and then filtering by the Bruhat strata in $\Bbf \backslash \Pbf w \Pbf /\Bbf$ yields the lemma.
\end{proof}

\section{Disjointness of the cohomology}

\subsection{Roots of the full twist}

\begin{defi}[\cite{SpringerRegularElements}]
An element $w \in \widetilde{\Wbf}$ is regular if it has a nonzero eigenvector in $V_{\C, \reg}$. An element is called $d$-regular if it is regular and has order $d$ and given $\zeta_d$ some primitive root of $1$, an element is called $\zeta_d$-regular if it has a regular eigenvector of eigenvalue $\zeta_d$. 
\end{defi}

\begin{thm}[\protect{\cite[Theorem 4.2, Theorem 6.4]{SpringerRegularElements}}]\leavevmode
\begin{enumerate}
\item All elements that are $d$-regular are $\zeta_d$-regular for any primitive $d$-th root of $1$.
\item There is at most one $\Wbf$-conjugacy class of $d$-regular elements in $\Wbf\delta \subset \widetilde{\Wbf}$. 
\item Let $w\delta$ be a $\zeta_d$-regular element, then $\Wbf(w\delta)$ is a complex reflection group with realization given by $\ker(w\delta - \zeta_d)$.
\end{enumerate}
\end{thm}

\begin{rque}
The values of $d$ for which there exist $d$-regular elements are classified in \cite[Section 6 for $\delta = 1$ and Section 8 for nontrivial $\delta$]{SpringerRegularElements}.
\end{rque} 

Let us fix some $d$ such that there exist $d$-regular elements in $\Wbf\delta$. We denote by $\widetilde{\Ocal}_d \subset \Wbf\cdot \delta$ the $\Wbf$-conjugacy class of $d$-regular elements.

\begin{defi}\label{defi:good-root-of-full-twist}
An element $c\delta \in \Wbf\delta$ is a root of the full twist if 
$$(\cbf\delta)^d = \FT \delta^d \in \widetilde{\Br}_{\Wbf}.$$
Moreover, we say that $c\delta$ is a \emph{good} root of the full twist if 
$$\ell((c\delta)^i) = i\ell(c\delta), \quad 1 \leq i \leq \left\lfloor \frac{d}{2} \right\rfloor.$$
We shall denote by 
$$\widetilde{\Ocal}_{d,\gd} \subset \widetilde{\Ocal}_d$$
the set of good roots of the full twist in $\widetilde{\Ocal}_d$.
\end{defi}

\begin{thm}[\protect{\cite[Théorème 3.12, Proposition 6.6]{BroueMichelActionRegularElements}}]
Let $c\delta$ be a root of the full twist, then $c\delta \in \widetilde{\Ocal}_d$.
\end{thm}

\begin{thm}[\protect{\cite[Théorème 3.17, Théorème 6.8]{BroueMichelActionRegularElements}}]\label{thm:all-good-roots-are-conjugate}
Any root of $\FT$ is conjugate in $\Br_{\Wbf}$ to a good root of $\FT$. 
\end{thm}

\begin{conj}[\protect{\cite[Conjecture 2.2]{DigneMichelEndomorphism}}]\label{conj:good-roots-of-full-twists-are-conjugate}
All good roots of $\FT$ in $\widetilde{\Ocal}_d$ are conjugate under $\Br_{\Wbf}$. 
\end{conj}

\begin{thm}[\protect{\cite[Theorem 12.4 (ii)]{Bessis}}]
Assume $\delta = 1$, then Conjecture \ref{conj:good-roots-of-full-twists-are-conjugate}.
\end{thm}

The goal of the rest of this section is to construct a specific element $c\delta \in \widetilde{\Ocal}_{d,\gd}$, i.e. one that satisfies Lemma \ref{lem:good-root-of-ft-selects-centralizer}. The main techniques for constructing such an element come from \cite{HeNie}.

\begin{lem}\label{lem:eigenvector-in-dominant-chamber}
There exists $c\delta \in \widetilde{\Ocal}_{d,\gd}$ such that $(c\delta + (c\delta)^{-1})$ has a nonzero eigenvector in $C^+$.
\end{lem}

\begin{proof}
We set $\zeta_d = \exp(\frac{2i\pi}{d})$. Let $c\delta \in \widetilde{\Ocal}_{d}$, since it is $\zeta_d$-regular, it has a nonzero eigenvector in $V_{\C, \reg}$. Let $v \in V_{\C, \reg}$ such that $(c\delta)v = \zeta_dv$. It follows that we have 
$$(c\delta + (c\delta)^{-1})v = 2\cos(2\pi/d)v$$
and therefore $\mathrm{Re}(v)$, the real part of $v$, is an eigenvector of $(c\delta + (c\delta)^{-1})$. Up to replacing $v$ by $e^{it}v$ for some $t \in [0, 2\pi]$, we can assume that $\mathrm{Re}(v) \in V_{\reg}$ and up to replacing $c$ by a conjugate, we can assume that $\mathrm{Re}(v)$ lies in $C^+$. 

We are now in the situation where $(c\delta + (c\delta)^{-1})$ has a nonzero eigenvector in $C^+$. Let $V_{\zeta_d} = \ker((c\delta + (c\delta)^{-1}) - 2\cos(2\pi/d), V) \subset V$. By assumption this contains a regular dominant element, it thus follows from \cite[Theorem 5.2]{HeNie}, applied to the element $c\delta$ and the sequence $r = (2\pi/d)$ that we have an equality in $\widetilde{\Br}_{\Wbf}$
$$(\cbf\delta)^d = \FT \delta^d.$$
Therefore the element $c\delta$ is a root of the full twist. It remains to see that this is a good root. By \cite[Lemma 2.1]{HeNie} applied to the element $(c\delta)^j$ for $j \leq  \left\lfloor \frac{d}{2} \right\rfloor$ and the subspace $V_{\zeta_d}$, we have 
$$\ell((c\delta)^j) = \frac{2j\pi}{d\pi}(|\Hfrak - \Hfrak_{V_{\zeta_d}}|) = j\frac{2\pi}{d\pi}(|\Hfrak - \Hfrak_{V_{\zeta_d}}|) = j\ell(c\delta),$$
and the element $c\delta$ is a good root of the full twist.
\end{proof}

\begin{lem}\label{lem:positivityDifference}
Let $u \leq v \in \Wbf$ and let $\lambda \in C^+$, then $(u - v)\lambda \in \sum_{\alpha \in \Delta} \Rbb_{\geq 0}\alpha$.
\end{lem} 

\begin{proof}
By the definition of the Bruhat ordering, see notation section \ref{ssec:bruhat-ordering}, there exists a sequence $(u_i)_{1 \leq i \leq n}$ and positive roots $(\alpha_i)$ such that 
\begin{enumerate}
\item $u_1 = u, u_n = v$ and for all $i$, we have $u_i < u_{i+1}$, 
\item $u_{i+1} = s_{\alpha_i}u_i$ and $u_i^{-1}\alpha_i > 0$. 
\end{enumerate}
We are going to prove that for all $i$, the difference $(u_i - u_{i+1})(\lambda)$ lies in $\sum_{\alpha \in \Delta} \Rbb_{\geq 0}\alpha$, taking the sum of all these terms yields the lemma. At a single step, we have 
\begin{align*}
(u_i - u_{i+1})\lambda &= (\id - s_{\alpha_i})(u_i\lambda) \\
&= 2\frac{\langle \alpha_i, u_i\lambda \rangle}{\langle \alpha_i, \alpha_i \rangle}\alpha_i \\
&= 2\frac{\langle u_i^{-1}\alpha_i, \lambda \rangle}{\langle \alpha_i, \alpha_i \rangle}\alpha_i.
\end{align*}
As $u^{-1}\alpha_i > 0$, we have $\langle u^{-1}\alpha_i, \lambda \rangle \geq 0$ and the element $(u_i - u_{i+1})\lambda$ lies in $\Rbb_{\geq 0}\alpha_i \subset \sum_{\alpha \in \Delta} \Rbb_{\geq 0}\alpha$.
\end{proof}

\begin{lem}\label{lem:good-root-of-ft-selects-centralizer}
Let $c\delta \in \Wbf\delta$ such that $c\delta + (c\delta)^{-1}$ has an eigenvector in $C^+$, then for all $x \in \Wbf$, if $xc \leq c\delta(x)$ and $\delta(x)c^{-1} \leq c^{-1}x$ then $xc = c\delta(x)$. 
\end{lem}

\begin{proof}
 Let $\rho \in C^+$ and $\alpha \in \Rbb$ be such that 
$$(c\delta + (c\delta)^{-1})(\rho) = \alpha\rho.$$

Let $x \in \Wbf$ such that $xc \leq c\delta(x)$ and $\delta(x)c^{-1} \leq c^{-1}x$. By Lemma \ref{lem:positivityDifference}, both 
$$(x(c\delta) - (c\delta)x)\rho \in \sum_{\alpha \in \Delta} \Rbb_{\geq 0}\alpha, \qquad (x(c\delta)^{-1} - (c\delta)^{-1}x)\rho \in \sum_{\alpha \in \Delta} \Rbb_{\geq 0}\alpha.$$
In particular, we have 
$$\left\langle (x(c\delta) - (c\delta)x)\rho, \rho \right\rangle \geq 0, \qquad \left\langle (x(c\delta)^{-1} - (c\delta)^{-1}x)\rho, \rho \right\rangle \geq 0.$$
Consider the sum 
\begin{align*}
&\left\langle (x(c\delta) - (c\delta)x)\rho, \rho \right\rangle + \left\langle (x(c\delta)^{-1} - (c\delta)^{-1}x)\rho, \rho \right\rangle  \\ 
&= \left\langle (x(c\delta) + x(c\delta)^{-1})\rho, \rho \right\rangle - \left\langle ((c\delta) + (c\delta)^{-1})x\rho, \rho \right\rangle \\
&= \alpha\langle x\rho, \rho\rangle - \alpha\langle \rho, x\rho \rangle \\
&= 0.
\end{align*}
As both terms are nonnegative, this forces $\langle (x(c\delta) - (c\delta)x)\rho, \rho \rangle = 0$ and as $\rho$ lies in the interior of the dominant chamber, this implies $(x(c\delta) - (c\delta)x)\rho = 0$ and as $\rho$ is regular, it has trivial stabilizer and thus $x(c\delta) = (c\delta)x$.
\end{proof}

\subsection{Roots of the relative full twist}

Let $\Pbf \subset \Gbf$ be a standard parabolic and let $I \subset S$ be the corresponding set of simple reflections. We set 
$$\FT_{\Wbf, I}  := \FT_{\Wbf} \FT_{\Wbf_{I}}^{-1} \in \Br_{\Wbf}.$$
This element is called the relative full twist. 

\begin{defi}[Periodic elements]\label{defi:periodic-elements}
A $(d,2k)$-periodic element is an element $b \in \Br_{\Wbf}^+\delta$ such that $I \subset S$ is stable under $\delta^d$ and
$$b^d = \delta^d\FT_{\Wbf,I}^k.$$
\end{defi} 

\begin{defi}[Non-extendable elements, \cite{DigneMichelParabolic}]
Let $b \in \Br_{\Wbf}^+\delta$ be a $(d,2k)$-periodic element, then $b$ is said to be non-extendable if there does not exist a $\delta^d$-stable proper subset $J \subsetneq I$ and $v \in \Br^+_{\Wbf_{I}}$ such that 
$$(vb)^d = \delta^d\FT_{\Wbf, J}^k.$$
\end{defi} 

\begin{conj}[\protect{\cite[Conjecture 8.9]{DigneMichelParabolic}}]\label{conj:non-extendable-elements-are-conjugate}
All non-extendable $(d,2k)$-periodic elements are conjugate.
\end{conj}

\begin{lem}\label{lem:existence-periodic-element}
Let $(d,k)$ be such that there exists a non-extendable $(d,2k)$-periodic element giving a root of $\FT_{\Wbf,I}$ with $2k\leq d$, then there exists $c \in \Wbf$ such that 
\begin{enumerate}
\item $\cbf\delta$ is a non-extendable $(d,2k)$-periodic element, 
\item $c\delta + (c\delta)^{-1}$ has an eigenvector in $\overline{C}^+_{I, \reg}$.
\end{enumerate}
\end{lem}

\begin{proof}
Let $\zeta = \exp(\frac{2ik\pi}{d})$ and let $c\delta$ be such that $\ker(c\delta - \zeta)$ is maximal. Let $v \in \ker(c\delta - \zeta)_{\reg}$ and set $\Wbf_I = \Stab_{\Wbf}(v)$. Arguing as in Lemma \ref{lem:eigenvector-in-dominant-chamber}, up to rescaling $v$ and conjugating $c$, we can assume that $\mathrm{Re}(v) \in \overline{C}^+_{I, \reg}$ and then replace $c$ by an element that is $I$-reduced on both sides. Moreover, $v$ is clearly an eigenvector of $c\delta + (c\delta)^{-1}$. It remains to see that $\cbf\delta$ is a non-extendable $(d,2k)$-periodic element. Note that by maximality of $\ker(c\delta - \zeta)$, by \cite[Theorem 8.3]{DigneMichelParabolic}, once we have established that $\cbf\delta$ is $(d,2k)$-periodic then it is non-extendable. We apply \cite[Lemma 5.1]{HeNie} to the element $c\delta$, let us check the hypotheses of the lemma: 
\begin{enumerate}
\item by the above construction, $\overline{C}^+$ contains a regular vector in $V_{\zeta} = \ker(c\delta + (c\delta)^{-1} - 2\cos(2k\pi/d)\id, V)$, 
\item we have $\Hfrak_{V_{\zeta}} = \{H_{\alpha}, \alpha \in \Phi_I^+\}$ and for all $v' \in V_{\zeta}$ and $\alpha \in \Phi_I^+$, we have 
$$\langle (c\delta)v, \alpha \rangle = \langle v, (c\delta)^{-1}\alpha \rangle \geq 0$$
and $c\delta$ is $I$-reduced on both sides $(c\delta)^{-1}(\Phi_I^+) \subset \Phi^+$. Consequently $C$ and $(c\delta)C$ are in the same connected component of $V_{\Rbb} \setminus \cup_{\alpha \in \Phi_I^+} H_{\alpha}$.
\end{enumerate}
By \cite[Lemma 5.1]{HeNie}, we thus have 
$$(\cbf\delta)^d = \delta^d \FT_{\Wbf,I}^k.$$
\end{proof}

\begin{lem}\label{lem:coset-centralizer-of-periodic-elements}
Let $c\delta \in \Wbf\delta$ be such that $\cbf\delta$ is a non-extendable $(d,2k)$-periodic element such that $(c\delta) + (c\delta)^{-1}$ has an eigenvector in $\overline{C}^+_{I, \reg}$ then for all $x \in \Wbf$ such that there exist $a,b \in \Wbf_{I}$ such that 
$$x(ac\delta) \leq (bc\delta)x, x(ac\delta)^{-1} \leq (bc\delta)^{-1}x,$$
then we have $x \in \Nbf_{\Wbf}(\Wbf_{I})$ and $xc\delta\Wbf_{I} = c\delta x\Wbf_{I}$. 
\end{lem}

\begin{proof}
Let $v \in \overline{C}^+_{I, \reg}$ be an eigenvector for $(c\delta) + (c\delta)^{-1}$. Note that as $c\delta$ normalizes $\Wbf_I$ we have $(ac\delta)v = (c\delta)(\Ad(c)\delta(a))v = (c\delta)v$ as $\Wbf_I$ acts trivially on $\overline{C}^+_{I, \reg}$. Similarly, we have $(bc\delta)v = (c\delta)v$. Arguing as in Lemma \ref{lem:good-root-of-ft-selects-centralizer}, we have 
$$\left\langle (x(ac\delta) - (bc\delta)x)v, v \right\rangle \geq 0$$
and 
$$\left\langle (x(ac\delta)^{-1} - (bc\delta)^{-1}x)v, v \right\rangle \geq 0.$$
As in Lemma \ref{lem:good-root-of-ft-selects-centralizer}, their sum is zero and therefore both terms are zero. In particular, we have 
$$\left\langle (x(c\delta) - (c\delta)x)v, v \right\rangle = 0.$$
Let $(y_i)_{i = 0, \dots, n}$ be a sequence of elements in $\Wbf\delta$ such that 
\begin{enumerate}
\item $y_0 = x(ac\delta)$, $y_n = (bc\delta)x$, 
\item $y_{i+1} > y_i$ and $y_{i+1} = s_{\beta_i}y_i$ where $\beta_i \in \Phi^+$ with $y_i^{-1}\beta_i > 0$.  
\end{enumerate}
Such a sequence exists by definition of the Bruhat ordering, see the notation Section \ref{ssec:bruhat-ordering}. Consider 
$$\left\langle (y_i - y_{i+1})v, v\right\rangle = 2\frac{\langle \beta_i, v\rangle\langle y_iv, \beta_i\rangle}{\langle \beta_i, \beta_i \rangle}.$$
By Lemma \ref{lem:positivityDifference} this is nonnegative and as 
$$\left\langle (x(ac\delta) - (bc\delta)x)v, v \right\rangle = \sum_i \left\langle (y_i - y_{i+1})v, v\right\rangle$$
all the terms must vanish. Therefore for all $i$, we have either $\langle \beta_i, v\rangle = 0$ or $\langle y_iv, \beta_i\rangle = 0$. 
\begin{enumerate}
\item if $\langle \beta_i, v\rangle = 0$, then $s_{\beta_i} \in \Wbf_I$, and $y_{i+1} \in \Wbf_Iy_i$, 
\item if $\langle y_iv, \beta_i\rangle = 0$, then $s_{y_i^{-1}\beta_i} \in \Wbf_I$ and $y_{i+1} \in y_{i}\Wbf_I$. 
\end{enumerate}
In both cases, we have $\Wbf_I y_i \Wbf_I = \Wbf_I y_{i+1} \Wbf_I$ and therefore 
\begin{equation}\label{eq:commutative-double-cosets}
\Wbf_I x(c\delta)\Wbf_I = \Wbf_I (c\delta)x \Wbf_I.
\end{equation}
Therefore, there exist $\alpha, \beta \in \Wbf_I$ such that $xc\delta = \alpha(c\delta)x\beta$. Let us rewrite this as 
$$x^{-1}\alpha c\delta x = c\delta \beta^{-1}.$$
Let $\zeta = \exp(\frac{2ik\pi}{d})$, by \cite[Theorem 8.3, (iii)]{DigneMichelParabolic}, as $c\delta$ is non-extendable, the $\zeta$-eigenspaces of $c\delta$ and of $\gamma c\delta$ coincide for any $\gamma \in \Wbf_I$. Consequently, $x$ normalizes $\ker(c\delta - \zeta)$ and by \cite[Lemma 8.7]{DigneMichelParabolic}, the normalizer of $\ker(c\delta - \zeta)$ is $\Nbf_{\Wbf}(\Wbf_I)$. Hence both $x$ and $c\delta$ normalize $\Wbf_I$ and by \eqref{eq:commutative-double-cosets}, we have 
$$x(c\delta)\Wbf_I = (c\delta)x \Wbf_I.$$
\end{proof}

\subsection{Mapping spaces between standard objects}

\begin{lem}\label{lem:filtration-product-of-Delta}
Let $x,c$ be two elements in $\Wbf$. The object $\Delta_x\Delta_c$ has a filtration in $\Hcal$ whose graded pieces are direct summands of $\Delta_y[i]$ for some $i \in \Z$ and $y \geq xc$. Moreover, $\Delta_{xc}[i]$ appears only for $i = 0$ and with multiplicity one. 
\end{lem}

\begin{proof}
This is a geometric reformulation of \cite[Lemma 8.9]{DigneMichelEndomorphism}. Let $y \in \Wbf$ such that $\Bbf y \Bbf \cap \Bbf x \Bbf c \Bbf \neq \emptyset$, then $\Bbf x^{-1}\Bbf \cap \Bbf c\Bbf y^{-1}\Bbf \neq \emptyset$ and by \cite[Chapitre IV, Section 2, Lemme 1]{BourbakiLie456}, we have 
$$\Bbf c\Bbf y^{-1}\Bbf  \subset \bigcup_{t \leq y^{-1}} \Bbf ct \Bbf.$$
There exists $t \leq y^{-1}$ such that $x^{-1} = ct$, i.e. $xc = t^{-1} \leq y$. 
As the convolution map factors 
$$\Ubf\backslash \Bbf x \Bbf /\Ubf \times \Ubf \backslash \Bbf c \Bbf /\Ubf \to \Ubf\backslash \Bbf x \Bbf c \Bbf /\Ubf \xrightarrow{j} \Ubf \backslash \Gbf/\Ubf$$
and any object in the essential image of $j_!$ has a filtration whose graded pieces are direct summands of $\Delta_y[i]$ with $\Bbf y \Bbf \cap \Bbf x \Bbf c \Bbf \neq \emptyset$, i.e. $y \geq xc$.

Let us now prove the equality case. Consider the Cartesian square 
\[\begin{tikzcd}
	\Xcal_{\Ubf} & {\Ubf\backslash\Bbf x\Bbf \times^{\Ubf} \Bbf c\Bbf/\Ubf} \\
	{\Ubf\backslash\Bbf xc\Bbf/\Ubf} & {\Ubf\backslash\Gbf/\Ubf}
	\arrow[from=1-1, to=1-2]
	\arrow[from=1-1, to=2-1]
	\arrow["m", from=1-2, to=2-2]
	\arrow[from=2-1, to=2-2]
\end{tikzcd}\]
Let $z = xc$ and let us now compute the fiber at the point $\dot{z} := \dot{x}\dot{c}$ of $\Delta_x\Delta_c$, we show that $\dot{z}^*\Delta_x\Delta_c \cong \oplus_{s\in\Ch_{\Lambda}(\Tbf)} (\Lcal \otimes_{\Lambda} \Lcal_s) [\dim(\Tbf) + \ell(xc)]$, since $\Ubf\backslash\Bbf xc\Bbf/\Ubf$ is closed in the support of $\Delta_x\Delta_c$, this implies that 
$$\Hom(\Delta_x * \Delta_c, \Delta_{xc}) \cong \prod_{s \in \Ch(\Tbf)}\RR[0],$$
and the filtration we constructed only has one copy of $\Delta_{xc}$ in degree $0$. Let us remove the left and right quotients by $\Ubf$ and consider the following diagram 
\[\begin{tikzcd}
	{\Xcal_{\Ubf, \dot{z}}} & {\Bbf x\Bbf \times^{\Ubf} \Bbf c \Bbf} \\
	{\Xcal_{\Bbf, \dot{z}}} & {\Bbf x\Bbf \times^{\Bbf} \Bbf c \Bbf} \\
	{\dot{z}} & \Gbf
	\arrow[from=1-1, to=1-2]
	\arrow["{m_{1,z}}"', from=1-1, to=2-1]
	\arrow["{m_1}", from=1-2, to=2-2]
	\arrow[from=2-1, to=2-2]
	\arrow["{m_{2,z}}"', from=2-1, to=3-1]
	\arrow["{m_2}", from=2-2, to=3-2]
	\arrow[from=3-1, to=3-2]
\end{tikzcd}\]
Introduce $\Ubf_x = \Ubf \cap x\Ubf x^{-1}$ and $\Ubf^{x} =  \Ubf \cap x\Ubf^- x^{-1}$ so that we have $\Bbf x\Bbf = \Ubf^x x\Bbf$ and the product is unique. Consequently, we have 
$$\Bbf x\Bbf \times^{\Bbf} \Bbf c \Bbf = \Ubf^x \times x\Bbf c \times \Ubf^{c^{-1}},$$
and 
\begin{align*}
&m_{1,!}(\mu_x^*(\bigoplus_{s} (\Lcal \otimes_{\Lambda} \Lcal_s)[\dim (\Tbf) + \ell(x)]) \boxtimes  \mu_c^*(\bigoplus_{s'} (\Lcal \otimes_{\Lambda} \Lcal_{s'})[\dim (\Tbf) + \ell(c)])) \\
&\cong \mu_{x,c}^*(\bigoplus_{s} (\Lcal \otimes_{\Lambda} \Lcal_s)[\dim (\Tbf) + \ell(x) + \ell(c)])
\end{align*}
where $\mu_x : x\Bbf \to \Tbf$ is the map induced by the choice of $\dot{x}$ and $\mu_{x,c} : x\Bbf c \to \Tbf$ is the map induced by the choice of $\dot{c}$ and $\dot{x}$.

Let us now compute the fiber $\Xcal_z := \Xcal_{\Bbf,\dot{z}}$. Writing $\Bbf = \Tbf\Ubf$, this fiber is the collection of tuples $(u^x, t, u, u^{c^{-1}}) \in \Ubf^x \times \Tbf \times \Ubf \times \Ubf^{c^{-1}}$ subject to the equation
\begin{equation}\label{eq:equality-case}
u^x \dot{x} t u \dot{c} u^{c^{-1}} = \dot{z}.
\end{equation} 
It follows immediately that \eqref{eq:equality-case} forces $t = 1$ and therefore, the fiber over $\dot{z}$ of $\Delta_x\Delta_c$ is $\oplus_{s\in\Ch_{\Lambda}(\Tbf)} (\Lcal \otimes_{\Lambda} \Lcal_s) [\dim(\Tbf) + \ell(x) + \ell(c)] \otimes \RGamma_c(\Xcal_z, \Lambda)$. 

We can rewrite \eqref{eq:equality-case} as 
$$\dot{x} \Ad(\dot{x}^{-1})(u^x) u \Ad(\dot{c})(u^{c^{-1}}) \dot{c} = \dot{z}$$
and therefore
\[\Ad(\dot{x}^{-1})(u^x) u \Ad(\dot{c})(u^{c^{-1}}) = 1.\]
Since both $\Ad(\dot{x}^{-1})(u^x)$ and $\Ad(\dot{c})(u^{c^{-1}})$ lie in $\Ubf^-$, we have $u = 1$ and 
$$\Ad(\dot{x}^{-1})(u^x) = \Ad(\dot{c})(u^{c^{-1}})^{-1}.$$
It then follows that we have an isomorphism 
$$\Xcal_z \cong \Ubf \cap \Ad(\dot{x})(\Ubf^-) \cap \Ad(\dot{z})(\Ubf).$$
By Lemma \ref{lem:lenght-intersection-of-unip-groups}, we have 
$$\dim(\Ubf \cap \Ad(\dot{x})(\Ubf^-) \cap \Ad(\dot{z})(\Ubf)) = (\ell(x) + \ell(c) - \ell(z))/2$$ 
and therefore $\RGamma_c(\Xcal_z, \Lambda) \cong \Lambda[-\ell(x) - \ell(c) + \ell(z)]$ and the shifts cancel out as expected.
\end{proof}

\begin{lem}\label{lem:lenght-intersection-of-unip-groups}
Let $x,y,z \in \Wbf$ with $z = xy$. Then we have 
$$\dim(\Ubf \cap \Ad(\dot{x})(\Ubf^-) \cap \Ad(\dot{z})(\Ubf)) = (\ell(x) + \ell(y) - \ell(z))/2.$$
\end{lem}

\begin{proof}
For signs $\epsilon, \delta \in \{\pm 1\}$, set 
$$\Phi_{\epsilon,\delta} := \{\alpha \in \Phi^+, x^{-1}\alpha \in \Phi^{\epsilon}, z^{-1}\alpha \in \Phi^{\delta}\}.$$
We have 
$$\dim(\Ubf \cap \Ad(\dot{x})(\Ubf^-) \cap \Ad(\dot{z})(\Ubf)) = |\Phi_{-,+}|.$$
Now observe that there is an equality 
$$\{\alpha \in \Phi^+, y^{-1}\alpha < 0 \} = x^{-1}\bigl((-\Phi_{-,+}) \bigsqcup \Phi_{+,-}\bigr)$$
and in particular we have $\ell(y) = |\Phi_{+,-}| + |\Phi_{-,+}|$. 
Finally, we have
\begin{align*}
\ell(x) + \ell(y) - \ell(z) &= (|\Phi_{-,+}| + |\Phi_{-,-}|) + (|\Phi_{-,+}| + |\Phi_{+,-}|) - (|\Phi_{+,-}| + |\Phi_{-,-}|) \\
&= 2|\Phi_{-,+}|.
\end{align*}
\end{proof}

\begin{rque}
Lemma \ref{lem:lenght-intersection-of-unip-groups} seems to be classical; however, we have not been able to find it as stated in the literature (even using AI to do so). 
\end{rque}

\begin{lem}\label{lem:filtration-product-of-delta-bis}
Let $x,y \in \Wbf$, then the product $\Delta_x\Delta_y$ has a filtration whose graded pieces are of the form $\Delta_z[j]$ with $z \in \Wbf$ and $j \leq 0$. Similarly $\nabla_x\nabla_y$ has a filtration whose graded pieces are summands of $\nabla_z[j]$ with $z \in \Wbf$ and $j \geq 0$.
\end{lem}

\begin{rque}
Note that in Lemma \ref{lem:filtration-product-of-delta-bis}, unlike in Lemma \ref{lem:filtration-product-of-Delta}, we do not impose $z \geq xy$, but we do impose $j \leq 0$. The argument for Lemma \ref{lem:filtration-product-of-delta-bis} is by induction and less geometric. 
\end{rque}

\begin{proof}
We argue by induction. Writing $y = s_1\dots s_n$ as a reduced expression, we can assume that $y = s_1 = s$ is a simple reflection. Then if $\ell(xs) = \ell(x) + 1$, we have $\Delta_x\Delta_s = \Delta_{xs}$. Otherwise, we have $\ell(xs) = \ell(x) - 1$ and thus $\Delta_x\Delta_s = \Delta_{xs}\Delta_s\Delta_s$, we are therefore reduced to the case $x = y = s$. By Lemma \ref{lem:product-simple-reflection}, we can filter $\Delta_s\Delta_s$ with graded pieces that are summands of $\Delta_s, \Delta_s[-1]$ and $\Delta_1$ hence the lemma is proved. The case of $\nabla_x\nabla_y$ is analogous.
\end{proof}

\begin{lem}\label{lem:filtration-product-of-delta-nabla}
Let $x,c$ be two elements in $\Wbf$. The object $\Delta_c\nabla_{x}$ has a filtration in $\Hcal$ whose graded pieces are direct summands of $\Delta_z[j]$ for $j \in \Z$ and $z \leq cx$.
Moreover, $\Delta_{cx}[j]$ appears only with $j = 0$ and multiplicity one.
\end{lem}

\begin{proof}
The element $b = \cbf\overline{\xbf}^{-1} \in \Br_{\Wbf}$, where $\overline{\xbf}$ is the braid lift of $x^{-1}$, is a rational permutation braid in the sense of \cite[Definition 4.2]{DigneGobet}. By \cite[Lemma 4.4]{DigneGobet}, given a reduced expression $cx = s_1\dots s_n$ in $\Wbf$, there exist signs $\epsilon_i \in \{\pm 1\}$ such that
$$b = \sbf_1^{\epsilon_1} \dots \sbf_n^{\epsilon_n}.$$
Set $w = \pi(b)$. We now prove that for any such braid $b$, the object $\Delta_b$ has a filtration whose graded pieces are of the form $\Delta_z[j]$ with $z \leq \pi(b)$.

First consider the case of a simple reflection $s$ and define
$$K_s := \mathrm{cone}(\Delta_s \to \nabla_s).$$
By Lemma \ref{lem:structure-simple-reflection}, the object $K_s$ has a filtration whose graded pieces are summands of $\Delta_1$ and $\Delta_1[1]$. Correspondingly $\nabla_s$ has a filtration whose graded pieces are $\Delta_s$ and $K_s$. 

Convolving the above triangle by $\Delta_s$ also shows that $\Delta_s * \Delta_s$ has a filtration whose graded pieces are direct summands of $\Delta_s[i]$ and $\Delta_1[i']$. 

We now proceed by induction on $n$. If $n = 0$, the statement is clear. Assume now that $n > 0$, set $b' = \sbf_1^{\epsilon_1} \dots \sbf_{n-1}^{\epsilon_{n-1}}$ and $w' = \pi(b')$ so that $w = w's_n$ and $b = b'\sbf_n^{\epsilon_n}$. If $\epsilon_n = -1$, by the case of a simple reflection, we can filter $\nabla_{s_n}$ and replace it by $\Delta_{s_n}$ and $K_{s_n}$ hence passing to graded we can replace $b$ by $b'\sbf_n$ and $b'$. We can therefore assume that $\epsilon_n = 1$. By induction, we can filter $\Delta_{b'}$ with graded pieces isomorphic to direct summands of $\Delta_{z}[j]$ with $z \leq w'$, multiplying by $\Delta_{s_n}$ on the right shows that $\Delta_b$ has a filtration whose graded pieces are direct summands of $\Delta_z * \Delta_{s_n}[j]$. Given such a $z$: 
\begin{enumerate}
\item if $\ell(zs_n) = \ell(z) + 1$, then $\Delta_z * \Delta_{s_n}[j] = \Delta_{zs_n}[j]$, 
\item if $\ell(zs_n) = \ell(z) -  1$, then $\Delta_z * \Delta_{s_n} = \Delta_{zs_n} * \Delta_{s_n} * \Delta_{s_n}$. Using the filtration of $\Delta_{s_n} * \Delta_{s_n}$ we can replace $\Delta_{s_n} * \Delta_{s_n}$ by direct summands of $\Delta_{s_n}$ and $\Delta_1$ and therefore $\Delta_z * \Delta_{s_n}$ has a filtration whose graded pieces are direct summands of $\Delta_{zs_n}$ and $\Delta_{z}$ (up to shifts), as both $z \leq w' \leq w$ and $zs_n \leq w$, the induction is complete.
\end{enumerate}

Let us now prove that $\Delta_{cx}$ appears only in degree $0$ and with multiplicity one. If $n = 0$, this is clear. Assume $n > 0$ and set $w', b'$ as above. By induction, $\Delta_{w'}$ appears as a graded piece of $\Delta_{b'}$ with multiplicity one. Multiplying by $\Delta_{s_n}$ shows that $\Delta_w$ appears as a graded piece in degree $0$. Let us now show that this is the only contribution. By construction, the other terms are obtained by filtering $\Delta_z *\Delta_{s_n}$ where $z < w'$, but those can contribute only copies of $\Delta_z$ and $\Delta_{zs_n}$ both of which are $< w$. 
\end{proof}

\begin{lem}\label{lem:control-mapping-space}
Let $x,c$ be two elements in $\Wbf$, assume that we have 
$$\Hom(\Delta_x\Delta_c, \Delta_c\nabla_{\delta(x)}) \neq 0$$
then we have $xc \leq c\delta(x)$ and $\delta(x)c^{-1} \leq c^{-1}x$. 
\end{lem}

\begin{proof}
By Lemma \ref{lem:filtration-product-of-Delta}, the object $\Delta_x\Delta_c$ has a filtration whose graded pieces are direct summands of $\Delta_y[i]$ with $y \geq xc$ and by Lemma \ref{lem:filtration-product-of-delta-nabla}, the object $\Delta_c\nabla_{\delta(x)}$ has a filtration whose graded pieces are direct summands of $\Delta_z[j]$ with $z \leq c\delta(x)$. 
It follows that $\Hom(\Delta_x\Delta_c, \Delta_c\nabla_{\delta(x)})$ has a filtration whose graded pieces are direct summands of $\Hom(\Delta_y,\Delta_z)$ (up to shifts) for $xc \leq y$ and $z \leq c\delta(x)$. If $\Hom(\Delta_x\Delta_c, \Delta_c\nabla_{\delta(x)})$ is nonzero, then at least one of the graded pieces is nonzero, consequently, there exist $y$ and $z$ such that $\Hom(\Delta_y, \Delta_z) \neq 0$. This implies $y \leq z$ and thus we have 
$$xc \leq y \leq z \leq c\delta(x).$$
Hence the first inequality is proved.

For the second inequality, observe that 
$$\Hom(\Delta_x\Delta_c, \Delta_c\nabla_{\delta(x)}) \cong \Hom(\nabla_{c^{-1}}\Delta_x, \nabla_{\delta(x)}\nabla_{c^{-1}}).$$
Dualizing Lemma \ref{lem:filtration-product-of-Delta}, we get that $\nabla_{\delta(x)}\nabla_{c^{-1}}$ has a filtration by direct summands of $\nabla_z[j]$ with $z \geq \delta(x)c^{-1}$ and by Lemma \ref{lem:filtration-product-of-delta-nabla}  the object $\nabla_{c^{-1}}\Delta_x$ has a filtration whose graded pieces are direct summands of $\Delta_y[i]$ with $y \leq c^{-1}x$. If the mapping space is nonzero, then for some $y,z$, we have 
$$\Hom(\Delta_y, \nabla_z) \neq 0$$
and in particular $\delta(x)c^{-1} \leq y = z \leq c^{-1}x$.
\end{proof}

\begin{lem}\label{lem:control-mapping-space-equality-case}
Let $x,c$ be two elements in $\Wbf$ and assume that $xc = c\delta(x)$, then 
$$\Hom(\Delta_x\Delta_c, \Delta_c\nabla_{\delta(x)}) \cong \prod_{s \in \Ch(\Tbf)}\RR_{\Lambda}[0].$$
\end{lem}

\begin{proof}
As in the proof of Lemma \ref{lem:control-mapping-space}, we can filter $\Hom(\Delta_x\Delta_c, \Delta_c\nabla_{\delta(x)})$ with graded pieces isomorphic to direct summands of $\Hom(\Delta_y, \Delta_z)[j]$ with $xc \leq y \leq z \leq c\delta(x)$. As $xc = c\delta(x)$, the only relevant pieces are $\Delta_{xc}$ in the filtration of $\Delta_x\Delta_c$ and $\Delta_{c\delta(x)}$ in the filtration of $\Delta_c\nabla_{\delta(x)}$. But by the equality cases of Lemmas \ref{lem:filtration-product-of-Delta} and \ref{lem:filtration-product-of-delta-nabla} these appear with multiplicity one in degree $0$.
\end{proof}

\subsection{Concentration in degree $0$: the Borel case}

\begin{thm}\label{thm:concentration-degree-0}
There exists $c\delta \in \widetilde{\Ocal}_{d,\gd}$, then 
$$\End(\RGamma_c(Y(\dot{c}), \Lambda))$$
lies in degree $0$ and in particular over $\Lambda = \Qlb$, 
$$\dim(\End(\RGamma_c(X(c), \Lambda))) = |\Wbf(c\delta)|.$$
\end{thm}

\begin{proof}
By Lemmas \ref{lem:eigenvector-in-dominant-chamber} and \ref{lem:good-root-of-ft-selects-centralizer}, there exists an element $c$ such that for all $x \in \Wbf$, if $xc \leq c\delta(x)$ and $\delta(x)c^{-1} \leq c^{-1}x$ then $xc = c\delta(x)$, let us choose this element. 

By adjunction, we have 
$$\End(\RGamma_c(Y(\dot{c}), \Lambda)) \cong \Hom(\hc_!\ch \pfrak_!\Delta_c, \pfrak_!\Delta_c)[-2\dim(\Ubf)]$$
and by Lemma \ref{lem:mackayfiltration}, $\hc^!\ch$ is filtered with graded pieces isomorphic to $\Int_x^!$. By Lemma \ref{lem:intertwinning-operators}, there is a natural isomorphism
\[\Int_x^! \pfrak_! \cong \pfrak_!(\Delta_x * - * \Delta_{\delta(x)^{-1}})[-2\dim(\Ubf)]\]
and in particular, the shifts cancel out. By Lemma \ref{lem:reduction-from-h-to-o}, we have 
$$\Hom(\pfrak_!(\Delta_x * \Delta_c * \Delta_{\delta(x)^{-1}}), \pfrak_! \Delta_c) \cong \Hom(\Delta_x * \Delta_c * \Delta_{\delta(x)^{-1}}, \Delta_c) \otimes_{\RR \otimes \RR} \RR.$$
By Lemmas \ref{lem:control-mapping-space} and \ref{lem:control-mapping-space-equality-case}, this mapping space is either $0$ if $xc \neq c\delta(x)$ or isomorphic to $\RR$ in degree $0$ otherwise. Moreover, as a bimodule it is isomorphic to $\RR_c$ and therefore $\RR_c \otimes_{\RR \otimes \RR} \RR_{\Frob}$ is in degree $0$. The case of $X(c)$ over $\Qlb$ is clear by Remark \ref{rque:series}.
\end{proof}

\begin{corol}
Over $\Lambda = \Qlb$, for all $i \neq j$ and $c\delta \in \widetilde{\Ocal}_{d,\gd}$, we have 
$$\left\langle H^i_c(Y(\dot{c}), \Qlb), H^j_c(Y(\dot{c}), \Qlb) \right\rangle = 0.$$
\end{corol}

\begin{proof}
Immediate from Theorem \ref{thm:concentration-degree-0}.
\end{proof}

\begin{thm}\label{thm:surjectivity}
Assume that $\delta = \id$. Let $c\delta \in \widetilde{\Ocal}_{d, \gd}$. The map 
$$\Qlb[\Br_{\Wbf}(\cbf\delta)] \to \End(\RGamma_c(X(c), \Qlb))$$
is surjective. 
\end{thm}

\begin{rque}
In the proof of Theorem \ref{thm:surjectivity}, we essentially prove that
\[\End(\RGamma_c(X(c), \Qlb))\]
is isomorphic to some Hecke algebra attached to the complex reflection group $\Wbf(c)$. However, we do not establish that this algebra is a cyclotomic Hecke algebra in the sense of \cite{BroueMalleCyclotomicHeckeAlg}. 
\end{rque}

\begin{proof}
The group $\Wbf(c)$ is a complex reflection group by \cite{SpringerRegularElements} and by \cite{Bessis}, there is a natural isomorphism 
$$\Br_{\Wbf(c)} \cong \Br_{\Wbf}(\cbf).$$
In particular, the map 
\begin{equation}\label{eq:braid-action}
\Qlb[\Br_{\Wbf}(\cbf)] \to \End(\RGamma_c(X(c), \Qlb))
\end{equation}
then defines a map from the group algebra of $\Br_{\Wbf(c)}$. We are going to prove that for any simple reflection $s \in \Wbf(c)$ of order $e$ the corresponding element $D_s \in \End(\RGamma_c(X(c), \Qlb))$ has $e$ distinct eigenvalues when it acts on $\RGamma_c(X(c), \Qlb)$ and spans a subalgebra of dimension $e$. These two facts imply that the minimal polynomial of $D_s$ in the endomorphism algebra is 
$$(X - \lambda_1)\dots (X - \lambda_e)$$
where $\lambda_1, \dots, \lambda_e$ are the distinct eigenvalues of $D_s$. It follows that the map \eqref{eq:braid-action} factors through a Hecke algebra of $\Wbf(c)$ with parameters $(\lambda_1, \dots, \lambda_e)$ for the reflection $s$, in the sense of \cite{BroueMalleCyclotomicHeckeAlg}. As these Hecke algebras have dimension $|\Wbf(c)|$ by \cite[4.C]{BroueMalleRouquier}, the surjectivity follows. 

Let $\Hfrak_d = \{H \cap \ker(c - \zeta_d), H \in \Hfrak \}$, by \cite{SpringerRegularElements}, this is the set of reflecting hyperplanes of $\Wbf(c)$. Given $K \in \Hfrak_d$, we denote by $\Wbf_{K} = \Stab_{\Wbf}(K)$ its pointwise stabilizer. This is a parabolic subgroup of $\Wbf$ and $\Wbf_K \cap \Wbf(c)$ is a rank-one reflection group. In \cite{DigneMichelEndomorphism}, in the last paragraph of section 3, the authors construct an element $t_K \in \Wbf_K \cap \Wbf(c)$ with the following properties 
\begin{enumerate}
\item $t_K$ is a generator of $\Wbf_K \cap \Wbf(c)$, 
\item $\tbf_K \in \Br_{\Wbf_K}$ satisfies $\tbf_K^e = \FT_{\Wbf_K}$.
\end{enumerate}
Let us denote by $\Lbf_K \subset \Pbf_K$ the standard Levi and parabolic subgroups attached to $\Wbf_K$. The group $\Wbf_K$ is normalized by $c$ and correspondingly, the group $\Lbf_K$ is normalized by $\dot{c}\Frob$. Let $Y(\dot{c}, \Pbf_K)$ be the corresponding parabolic Deligne--Lusztig variety. By transitivity of Deligne--Lusztig induction, we have 
$$Y(\dot{c}) \cong Y(\dot{c}, \Pbf_K) \times^{\Lbf^{\dot{c}\Frob}} Y^{\Lbf_K}(\dot{c})$$ 
where $Y^{\Lbf_K}(\dot{c})$ is the Deligne--Lusztig variety of $\Lbf_K$ attached to the element $\dot{c}$. Moreover, by \cite[Proposition 5.2.8]{DigneMichelRouquier}, this isomorphism is compatible with the action of the operators $D_b$, i.e. we have a commutative diagram 
\[\begin{tikzcd}
	{\Qlb[\Br_{\Wbf_K}(\cbf)]} & {\Qlb[\Br_{\Wbf}(\cbf)]} \\
	{\End(\RGamma_c(Y^{\Lbf_K}(\dot{c}), \Qlb))} & {\End(\RGamma_c(Y(\dot{c}), \Qlb)).}
	\arrow[from=1-1, to=1-2]
	\arrow[from=1-1, to=2-1]
	\arrow[from=1-2, to=2-2]
	\arrow[from=2-1, to=2-2]
\end{tikzcd}\]
By Theorem \ref{thm:concentration-degree-0} applied to the group $\Lbf_K$, $\End(\RGamma_c(Y^{\Lbf_K}(\dot{c}), \Qlb))$ is in degree $0$ and of dimension $e$. By \cite[Proposition 5.2.3]{DigneMichelRouquier}, for all $1 \leq i < e$, the operator $D_{\tbf_K^i}$ has no fixed points on $Y^{\Lbf_K}(\dot{c})$ and by \cite[Proposition 5.2.2]{DigneMichelRouquier}, we have 
$$\tr(D_{\tbf_K^i}, \RGamma_c(Y^{\Lbf_K}(\dot{c}), \Qlb)) = 0.$$
This implies that the operator $D_{\tbf_K}$ has at least $e$ distinct eigenvalues when acting on $\RGamma_c(Y^{\Lbf_K}(\dot{c}), \Qlb)$ and the claim follows. 
\end{proof}

\subsection{Concentration in degree $0$: the parabolic case}

We fix the following data: 
\begin{enumerate}
\item $\Pbf = \Lbf\UP \subset \Gbf$ a standard parabolic subgroup and we denote by $I \subset S$ the corresponding subset of simple reflections, 
\item $c \in {^{\Pbf}}\Wbf^{\Frob(\Pbf)}$ an element such that $\Ad(\dot{c})\Frob(\Lbf) = \Lbf$ and $\cbf\delta$ is a non-extendable $(d,2k)$-periodic element which is a root of $\FT_{\Wbf, I}$ and such that $c\delta + (c\delta)^{-1}$ has an eigenvector in $\overline{C}^+_{I,\reg}$.
\end{enumerate}

Let us denote by $\Nbf_{I} \subset \Nbf_{\Wbf}(\Wbf_I)$ the subset of elements such that $n(\Phi_I^+) = \Phi_I^+$. We have 
$$\Nbf_{\Wbf}(\Wbf_{I}) = \Wbf_{I} \rtimes \Nbf_I.$$

\begin{lem}\label{lem:self-duality}
The algebra $\End_{\Gbf^{\Frob}}(\RGamma_c(Y(\dot{c}, \Pbf), \Lambda))$ is canonically self-dual.
\end{lem}

\begin{proof}
Let $R_{\dot{c}, \Pbf} = \RGamma_c(Y(\dot{c}, \Pbf), \Lambda)$. This is a perfect complex of $\Gbf^{\Frob}$-modules, so 
$$\End_{\Gbf^{\Frob}}(\RGamma_c(Y(\dot{c}, \Pbf), \Lambda)) \cong (R_{\dot{c}, \Pbf}^{\vee} \otimes_{\Lambda} R_{\dot{c}, \Pbf})^{\Gbf^{\Frob}}.$$
As $R_{\dot{c}, \Pbf}^{\vee} \otimes_{\Lambda} R_{\dot{c}, \Pbf}$ is a self-dual perfect complex of $\Gbf^{\Frob}$-representations, we have 
\begin{align*}
\End_{\Gbf^{\Frob}}(\RGamma_c(Y(\dot{c}, \Pbf), \Lambda))^{\vee} &\cong (R_{\dot{c}, \Pbf}^{\vee} \otimes_{\Lambda} R_{\dot{c}, \Pbf})^{\Gbf^{\Frob}, \vee} \\
&\cong (R_{\dot{c}, \Pbf}^{\vee} \otimes_{\Lambda} R_{\dot{c}, \Pbf})_{\Gbf^{\Frob}} \\
&\cong (R_{\dot{c}, \Pbf}^{\vee} \otimes_{\Lambda} R_{\dot{c}, \Pbf})^{\Gbf^{\Frob}},
\end{align*}
where the last isomorphism is provided by the norm map from coinvariants to invariants which is an isomorphism on perfect complexes.
\end{proof}

\begin{thm}\label{thm:vanishing-non-extendable-elements}
The endomorphism algebra 
$$\End_{\Gbf^{\Frob}}(\RGamma_c(Y(\dot{c}, \Pbf), \Lambda))$$ 
is in degree $0$. 
\end{thm}

\begin{rque}
Under Conjecture \ref{conj:non-extendable-elements-are-conjugate}, Theorem \ref{thm:vanishing-non-extendable-elements} holds for all non-extendable $(d,2k)$-periodic elements by \cite[Lemma 7.22]{DigneMichelParabolic}.
\end{rque}

\begin{proof}
By Lemma \ref{lem:comparison-with-dl-parabolic}, we want to prove that the endomorphism algebra
$$\End(\ch_{\Pbf}j_{c,!}^{\Pbf}(\Reg_{\Lbf^{\dot{c}\Frob}}))$$
is in degree $0$. By Lemma \ref{lem:mackay-filtration-parabolic-case}, we can filter this complex and the graded pieces are 
$$\Hom(\Int_{\Pbf, \Pbf,w}^!j_{c,!}^{\Pbf}(\Reg_{\Lbf^{\dot{c}\Frob}}), j_{c,!}^{\Pbf}(\Reg_{\Lbf^{\dot{c}\Frob}}))[2\dim(\UP)],$$
with $w \in {^\Pbf}\Wbf^{\Pbf}$.

We first prove that if $\Hom(\Int_{\Pbf, \Pbf,w}^!j_{c,!}^{\Pbf}(-), j_{c,!}^{\Pbf}(-))$, the mapping space in the category of functors, is nonzero then $w \in \Nbf_I$ and $w(c\delta) = (c\delta)w$. If this mapping complex is nonzero, then there exist $a,b \in \Wbf_{\Frob(\Pbf)}$ such that 
$$\Hom(\Int_{\Pbf, \Pbf,w}^!(\ch^{\Pbf}(\pfrak_!\Delta_c\Delta_a)), j_{c,!}^{\Pbf}(\pfrak_!\Delta_c\Delta_b)) \neq 0.$$
By Lemma \ref{lem:comparison-mackay}, we may filter this complex with graded pieces of the form (up to a shift)
$$\Hom(\Int_{x}^!(\pfrak_!(\Delta_c\Delta_a)), \pfrak_!(\Delta_c\Delta_b)).$$
By Lemma \ref{lem:intertwinning-operators}, we can rewrite this as
$$\Hom(\pfrak_!(\Delta_x\Delta_c\Delta_a\Delta_{\delta(x^{-1})}), \pfrak_!(\Delta_c\Delta_b)).$$
By Lemma \ref{lem:reduction-from-h-to-o}, if this mapping space is nonzero, then 
$$\Hom(\Delta_x\Delta_c\Delta_a, \Delta_c\Delta_b\nabla_{\delta(x)})$$
is nonzero. As we have $\ell(c) + \ell(a) = \ell(ca)$ and $\ell(c) + \ell(b) = \ell(cb)$, we have 
$$\Hom(\Delta_x\Delta_{ca}, \Delta_{cb}\nabla_{\delta(x)}) \neq 0.$$
Arguing as in Lemma \ref{lem:control-mapping-space}, we deduce that 
$$xca \leq cb\delta(x), \delta(x)(ca)^{-1} \leq (cb)^{-1}x.$$
By Lemma \ref{lem:coset-centralizer-of-periodic-elements}, we have $x \in \Nbf_{\Wbf}(\Wbf_I)$ hence the same is true for $w$ and we have $x(c\delta)\Wbf_{I} = (c\delta)x\Wbf_{I}$. As $w$ is the minimal coset representative, we have $w(c\delta) = (c\delta)w$. It follows that the only nonzero terms are those corresponding to $w \in \Nbf_{I}(c\delta)$.

To prove the concentration in degree $0$, by Lemma \ref{lem:self-duality}, it is enough to see that for all $X \in \DD^{\geq 0}(\Rep_{\Lambda}\Lbf^{\dot{c}\Frob})$ and $Y \in \DD^{\leq 0}(\Rep_{\Lambda}\Lbf^{\dot{c}\Frob})$, the complex 
$$\Hom(\ch_{\Pbf}j_{c,!}^{\Pbf}(X), \ch_{\Pbf}j_{c,!}^{\Pbf}(Y))$$
is concentrated in nonpositive degrees. Taking for $X$ and $Y$ the cohomology complexes of Deligne--Lusztig varieties (with and without compact support), which by \cite[Theorem 3.2.4]{EteveTilting} and \cite{ZhuCategoricalTame}, generate the positive and negative parts of $\DD(\Rep_{\Lambda}\Lbf^{\dot{c}\Frob})$ respectively, it is enough to prove that for all $a,b \in \Wbf_{I}$, the complex 
$$\Hom(\ch(\pfrak_!(\Delta_a\Delta_c)), \ch(\pfrak_!(\nabla_b\Delta_c)))$$
is concentrated in nonpositive degrees. Using Lemma \ref{lem:mackayfiltration} and Lemma \ref{lem:reduction-from-h-to-o}, we prove that 
$$\Hom(\Delta_w\Delta_a\Delta_c, \nabla_b\Delta_c\nabla_{\delta(w)})$$
has a filtration by copies of $\RR_{\Lambda}[j]$ with $j \geq 0$ so that after taking $\otimes_{\RR_{\Lambda} \otimes \RR_{\Lambda}} \RR_{\Lambda}$ we obtain a complex concentrated in nonpositive degrees. Note that by the preceding discussion, we can assume that $w \in \Nbf_{\Wbf}(\Wbf_I)$ and $wc\delta\Wbf_{I} = (c\delta)w\Wbf_I$. Let us write $w = \alpha n$ with $\ell(\alpha) + \ell(n) = \ell(w)$ and $\alpha \in \Wbf_I, n \in \Nbf_I(c\delta)$. Then we have 
$$\Delta_w\Delta_a\Delta_c = \Delta_{\alpha}\Delta_n\Delta_a\Delta_c = \Delta_{\alpha}\Delta_{nan^{-1}}\Delta_n\Delta_c$$
as $n$ normalizes $\Wbf_I$ and preserves the length in $\Wbf_I$. By the same argument, we also have  
$$\nabla_b\Delta_c\nabla_{\delta(w)} = \nabla_b\Delta_c\nabla_{\delta(\alpha)}\nabla_{\delta(n)} = \Delta_c\nabla_{c^{-1}bc}\nabla_{\delta(\alpha)}\nabla_{\delta(n)}.$$
By Lemma \ref{lem:filtration-product-of-delta-bis}, we can filter both $\Delta_{\alpha}\Delta_{nan^{-1}}$ and $\nabla_{c^{-1}bc}\nabla_{\delta(\alpha)}$ with graded pieces of the form $\Delta_y[i]$ for some $y \in \Wbf_{I}$ and $i \leq 0$ and $\nabla_z[j]$ with $z \in \Wbf_{\delta(I)}$ and $j \geq 0$ respectively. It follows that $\Hom(\Delta_w\Delta_a\Delta_c, \nabla_b\Delta_c\nabla_{\delta(w)})$ has a filtration whose graded pieces are 
$$\Hom(\Delta_y\Delta_n\Delta_c, \Delta_c\nabla_z\nabla_{\delta(n)})[j-i].$$
Assume this is nonzero and let us now rewrite 
$$\Delta_y\Delta_n\Delta_c = \Delta_{yn}\Delta_c$$
as $n$ is the minimal $\Wbf_I$-coset representative and 
$$\Delta_c\nabla_z\nabla_{\delta(n)} = \Delta_c\nabla_{z\delta(n)}.$$
Arguing as in Lemma \ref{lem:control-mapping-space}, we have 
$$ync \leq cz\delta(n), z\delta(n)c^{-1} \leq c^{-1}yn.$$
As we have $nc = c\delta(n) \in \Nbf_I$, these two inequalities imply respectively 
$$y \leq czc^{-1}, z \leq c^{-1}yc$$ 
and as $\Ad(c) : \Wbf_{\delta(I)} \to \Wbf_{I}$ is length preserving, we have 
$$y = czc^{-1}$$
and therefore
$$ync = cz\delta(n).$$
Arguing as in Lemma \ref{lem:control-mapping-space-equality-case}, we have 
$$\Hom(\Delta_y\Delta_n\Delta_c, \Delta_c\nabla_z\nabla_{\delta(n)}) = \RR_{\Lambda}[0].$$
It then follows that $\Hom(\Delta_w\Delta_a\Delta_c, \nabla_b\Delta_c\nabla_{\delta(w)})$ is concentrated in nonpositive degrees and the theorem is proved.
\end{proof}

\bibliographystyle{alpha}
\bibliography{biblio.bib}

\end{document}